\documentclass[10pt,a4paper]{amsart}

\usepackage{amssymb,stmaryrd,booktabs,euscript,enumerate}

\usepackage[pagebackref]{hyperref}
\usepackage{accents}

\makeatletter
\def\widebar{\accentset{{\cc@style\underline{\mskip10mu}}}}
\def\Widebar{\accentset{{\cc@style\underline{\mskip8mu}}}}
\makeatother

\theoremstyle{plain}
\newtheorem{theorem}{Theorem}[section]
\newtheorem{lemma}[theorem]{Lemma}
\newtheorem{proposition}[theorem]{Proposition}
\newtheorem{corollary}[theorem]{Corollary}

\newtheorem*{conjecture}{Conjecture}

\theoremstyle{definition}

\theoremstyle{remark}
\newtheorem{remark}[theorem]{Remark}

\newtheorem{example}[theorem]{Example}
\newtheorem*{ack}{Acknowledgments}

\newcommand{\Lie}[1]{\operatorname{\textsl{#1}}}
\newcommand{\lie}[1]{\operatorname{\mathfrak{#1}}}
\newcommand{\End}{\Lie{End}}

\newcommand{\Hol}{\Lie{Hol}}

\newcommand{\SL}{\Lie{SL}} 

\newcommand{\SO}{\Lie{SO}}
\newcommand{\so}{\lie{so}}
\newcommand{\SP}{\Lie{Sp}}

\newcommand{\SU}{\Lie{SU}}
\newcommand{\su}{\lie{su}}
\newcommand{\Un}{\Lie{U}}
\newcommand{\un}{\lie{u}}
\newcommand{\Gtwo}{\ifmmode{{\rm G}_2}\else{${\rm G}_2$}\fi}

\newcommand{\LC}{{\nabla^g}}
\newcommand{\RC}{\textup{R}^g} 
\newcommand{\Ric}{\textup{Ric}}
\newcommand{\Nt}{\widebar\nabla}

\newcommand{\Rt}{\widebar{\textup{R}}}
\newcommand{\bian}{\textup{b}}

\DeclareMathOperator{\tr}{tr}
\DeclareMathOperator{\coker}{coker}
\newcommand{\Hodge}{\mathord{\mkern1mu *}}
\newcommand{\hook}{\mathbin{\lrcorner}}
\newcommand{\bw}{\mathbin{\wedge}}
\newcommand{\abs}[1]{\left\lvert #1\right\rvert}
\newcommand{\norm}[1]{\left\lVert #1\right\rVert}
\newcommand{\real}[1]{\left\llbracket #1 \right\rrbracket}

\newcommand{\inp}[2]{\left< #1, #2\right>}
\newcommand{\vol}{\textup{vol}}

\newcommand{\T}{\checkmark}
\newcommand{\lb}[1]{\raisebox{-10pt}[12pt][0pt]{#1}}

\newcommand{\mW}{\EuScript{W}}
\newcommand{\mK}{\EuScript{K}}
\newcommand{\mR}{\EuScript{R}}
\newcommand{\mS}{\EuScript{S}}

\date{\today}
\begin{document}

\title[Curvature decomposition of $\Gtwo$ manifolds]%
{Curvature decomposition of $\mathbf{G_2}$ manifolds}
\date{\today}

\author{Richard Cleyton}
\address[Cleyton]{Humboldt-Universit\"at zu Berlin\\
  Institut f\"ur Mathematik\\
  Unter den Linden 6\\
  D-10099 Berlin\\
  Germany} 

\email{cleyton@mathematik.hu-berlin.de}

\author{Stefan Ivanov}
\address[Ivanov]{University of Sofia "St. Kl. Ohridski"\\
  Faculty of Mathematics and Informatics,\\
  Blvd. James Bourchier 5,\\
  1164 Sofia, Bulgaria,}
\address{and Max-Planck-Institut f\"ur Mathematik\\Vivatsgasse
  7\\D-53111 Bonn\\Germany} 
  
\email{ivanovsp@fmi.uni-sofia.bg}

\begin{abstract}
  Explicit formulas for the $G_2$-components of the Riemannian
  curvature tensor on a manifold with a $G_2$ structure are given in
  terms of Ricci contractions. We define a conformally invariant
  Ricci-type tensor that determines the $27$-dimensional part of the
  Weyl tensor and show that its vanishing on compact $G_2$ manifold with
  closed fundamental form forces the three-form to be parallel. A
  topological obstruction for the existence of a $G_2$ structure with
  closed fundamental form is obtained in terms of the integral norms
  of the curvature components.  We produce integral inequalities for
  closed $G_2$ manifold and investigate limiting cases. We make a study
  of warped products and cohomogeneity-one $G_2$ manifolds. As a
  consequence every Fern´andez-Gray type of $G_2$ structure whose scalar
  curvature vanishes may be realized such that the metric has holonomy
  contained in $G_2$.
\end{abstract}

\maketitle

\section{Introduction}
A $7$-dimensional manifold $M$ is called a $G_2$ manifold if its
structure group reduces to the exceptional Lie group $G_2$, i.e. in
the presence of a $G_2$ structure on $M$.  In practical terms this is
a non-degenerate three-form which is positive, in the sense that it
associates a definite metric to $M$.  This three-form is called the
\emph{fundamental form} of the $G_2$ manifold and is throughout this
paper denoted by $\phi$.  From the purely topological point of view, a
$7$-dimensional paracompact manifold is a $G_2$ manifold if and only
if it is an oriented spin manifold~\cite{LM}.  The reduction of the
structure group lifts the representation theory of $G_2$ to tensor
bundles and allows to split bundles and their sections according to
their $G_2$ decompositions.  It also equips $M$ with a unique
\emph{canonical $G_2$ connection} of minimal torsion, in addition to
the Levi-Civita connection of the associated metric.  The difference
of these is the \emph{intrinsic torsion} of the $G_2$ structure.
Decomposing the intrinsic torsion into $G_2$ irreducible parts
divides $G_2$ structures into $16$ torsion classes labeled by the
non-vanishing torsion components~\cite{FG}.  These components of
torsion are precisely the first order diffeomorphism invariants of the
$G_2$ structure $\phi$.
Based on the general theory of Cartan, Robert Bryant observes in
\cite{Br1} that, for a $G_2$ structure, the diffeomorphism invariants
which are polynomial in the derivatives of $\phi$ up to second order,
are sections of a vector bundle of rank~$392$.  In the same paper
Bryant describes the $G_2$ invariant splitting of this bundle in to
eleven $G_2$ irreducible components. This includes the splitting of
the Riemannian curvature tensor in to five components.  Three of these
are part of the Weyl tensor $W=W_{77}\oplus W_{64}\oplus W_{27}$
(subscripts denote dimensions of the corresponding bundle) and the
remaining two correspond to the traceless Ricci curvature and
scalar curvature of the underlying metric.

In the present note the focus is on the following basic question: What
type of conclusion about the intrinsic torsion may be obtained by
imposing conditions on the five curvature components of the Riemannian
curvature and vice-versa?

For instance, it is well known that when the torsion vanishes
completely, the Riemannian curvature is precisely
$W_{77}$,~\cite{MR0231313}.  This has several nice corollaries.
First, the observation by Bonan~\cite{Bo} that metrics with holonomy
$G_2$ are Ricci-flat follows immediately.  Secondly, the other four
curvature components are determined linearly by the covariant
derivative of the intrinsic torsion and a term of second order in the
components of the torsion, see Proposition~\ref{thm:curv}.

Moreover, the representation theory of $G_2$ determines precisely how
the square of the torsion and its derivatives contribute to curvature
components.  Qualitative results deduced from this are listed in
Tables~\ref{tab:g2-2} and~\ref{tab:g2-1}. On the quantitative level,
universal formulas for the scalar curvature and Ricci curvature in
terms of torsion and its covariant derivative must exist.  These were
computed for the special instance of $G_2$ structures with skew
torsion in~\cite{SITF} and for the general case in~\cite{Br1}.  For
obtaining such expressions as well as for their applicability the
following observations play crucial roles.
The fundamental form is parallel with respect to the Levi-Civita
connection precisely when the intrinsic torsion is zero.  If so, the
associated metric has holonomy contained in $G_2$ and we say that the
$G_2$ structure is \emph{parallel}.  Gray~\cite{Gr} showed (see
also~\cite{FG,Br,Salamon}) that a $G_2$ manifold is parallel precisely
when the fundamental form is harmonic.  This has the corollary that,
for a general $G_2$ structure, the intrinsic torsion is completely
determined by the differential and co-differential of the fundamental
three-form.  The relationship may be made explicit, see
section~\ref{sec:deriv-fund-three}, and so embeds the intrinsic
torsion in the exterior algebra.
The curvature component $W_{27}$ manifests itself in the form of a
second Ricci-type tensor on $G_2$-manifolds.  This has an expression
as a linear combination of the Ricci tensor and the $\star$-Ricci
tensor introduced in \cite{CI} (here re-baptized ``the $\phi$-Ricci
tensor'' to emphasize the dependency on the three-form $\phi$ rather
than the metric).  Using this and the fact that the space of symmetric
two-tensors on $\mathbb R^7$ injects into the space of three-forms
provides expressions for scalar curvature, Ricci curvature and
$W_{27}$, quadratic in the components of the derivative and
coderivative of the fundamental form and linear in exterior
derivatives (rather than covariant derivatives) of the same
components, see Lemma~\ref{prop:1}.  This is a generalization of the
equations previously found by Bryant~\cite{Br1} and
Friedrich-Ivanov~\cite{SITF}.  An equivalent formula involving the
canonical connection is also given, see Lemma~\ref{prop:new1}.  In
Theorem~\ref{thm:curv-decomp} we give the precise relations between
$W_{27}$ and the induced symmetric two-tensor, plus other algebraic
relations between curvature components, scalar curvature $s_g$ and the
Ricci tensor of the associated metric $g$.

Turning to applications we first re-examine the $G_2$ structures
$\phi$ is either parallel (zero torsion), nearly parallel (a
one-dimensional torsion component is non-zero) or locally conformally
equivalent to one of these types (allowing an additional non-zero
one-form part in the intrinsic torsion).  For nearly parallel
structures we recover results of Gray and Reyes-Carrion, i.e. that
the induced metric is Einstein with positive scalar
curvature~\cite{Gr2} and $W_{64} = 0 = W_{27}$~\cite{carrion}.
Conformal invariance of the Weyl curvature shows that the latter facts
hold for the wider class of $G_2$ manifolds locally conformal to
parallel and nearly parallel, see section~\ref{sec:1+4}.

In section~\ref{sec:2} we concentrate on $G_2$ manifolds with closed
fundamental three-form.  Closed $G_2$ structures are interesting for
several reasons.  In favorable circumstances a closed $G_2$ structure
on a compact manifold may be deformed to obtain a metric with holonomy
$G_2$.  This lies at the heart of the work on the existence of compact
parallel $G_2$ manifolds by Joyce~\cite{Joyce:book,JI&II} and
Kovalev~\cite{Kov}.  Furthermore, closed $G_2$ structures share
several features with nearly parallel $G_2$ structures.  Notably, in
both cases the full space of second order diffeomorphism invariants is
determined by the components of the Riemannian curvature of the
underlying metric, something that is, in general, not true for a $G_2$
structure, see~\cite{Br1} and Remark~\ref{2order}.

In~\cite{CI}, we found that a compact $G_2$ manifold with a closed
$G_2$ structure is Einstein only if it is parallel.  Using our new
formulas for the Ricci-type tensors this result is easily recovered,
allowing to interpret the `natural equations for closed $G_2$
structures' coined in~\cite{Br1}.  In particular we show that the
concept of closed $G_2$ structures with extremally pinched Ricci
curvature is precisely that of a closed $G_2$ for which the
$27$-dimensional component of the covariant derivative of the
intrinsic torsion is zero.

As mentioned above, general principles show that the curvature
component $W_{64}$ also has an expression in terms of torsion and its
derivatives, even for a general $G_2$ structure. This is however
hampered by the lack of a sensible representation of the constituent
tensors.  Yet, for closed $G_2$ structures an explicit formula is
available, see Lemma~\ref{prop:4}.  Certain constraints on the
curvature of a closed $G_2$ structure can only be satisfied when
$\phi$ is also coclosed.  To be precise, we have
\begin{theorem}\label{main1}
  Suppose $M$ is a compact manifold with closed $G_2$ structure
  $\phi$.  If $W_{27}=0$ then $\phi$ is parallel.  The same conclusion
  holds if $W_{77}$ and $W_{64}$ vanish simultaneously, in which case
  the associated metric is furthermore flat.
\end{theorem}
The first statement here is implicitly contained in~\cite{Br1}
although more is true.  On a compact $M$, the same conclusion
holds for the vanishing of any linear combination of the two
Ricci-type tensors, apart precisely from the case corresponding to an
extremally pinched Ricci metric.

From the structure equations of a closed $G_2$ structure it follows
that the intrinsic torsion is parallel with respect to the canonical
connection precisely when the associated metric is extremally Ricci
pinched and $W_{64}$ is zero.  There is a then a natural hierarchy of
closed $G_2$ structures. Those with zero torsion clearly form the
smallest class, contained in that of parallel torsion, which is in
turn a subclass of extremally Ricci pinched structures.

Compact $M$'s with closed $\phi$ have a natural diffeomorphism
invariant number $n(M,\phi)$ associated to them by taking the cup
product of the first Pontrjagin class $p_1(M)$ with the de Rham class
$[\phi]$ and pairing with the fundamental class $[M]$, i.e.
$n(M,\phi):=\langle p_1(M)\cup[\phi],[M]\rangle$.  Using Chern-Weil
theory and the decomposition of the Riemannian curvature we obtain an
expression for $n(M,\phi)$ through an integral of norms of curvature
components, see Theorem~\ref{thm:pontrj}.  Now each of the classes in
the hierarchy of closed $G_2$ structures are characterized by
realizing different equalities in a series of estimates for
$n(M,\phi)$ in terms of $L^2$ norms of curvature components.  For the
penultimate class we shall furthermore prove that such compact
manifold are locally isometric to the example constructed by Bryant
in~\cite{Br1}:
\begin{theorem}\label{thm:par-tors}
  Let $M$ be a $G_2$ manifold with a closed fundamental form $\phi$.
\begin{enumerate}[$(a)$]
\item Suppose the intrinsic torsion
  is parallel with respect to the canonical $G_2$ connection.  Then
  $(M,\phi)$ is locally isometric to the homogeneous space $G/\SU(2)$
  where $G$ is the Lie group consisting of affine transformations of
  two-dimensional complex space preserving its standard complex volume
  form.\label{item:4}
\item Suppose furthermore that $M$ is compact.  Then
  \begin{equation*}
    n(M,\phi) \geqslant -\frac1{8\pi^2}\int_M
    \left\{\norm{W_{77}}^2 - \tfrac{3}{16} s_g^2\right\} dV_g, 
  \end{equation*}
  and equality holds if and only if $M$ is as described in~(\ref{item:4}).
\end{enumerate}
\end{theorem}
In the light of the above results we make the following
\begin{conjecture}
  Suppose $M$ is compact and $\phi$ is a closed $G_2$ structure on $M$
  with extremally pinched Ricci curvature.  Then the universal
  covering space $\tilde M$ of $M$ is isometric to $G/\SU(2)$ equipped
  with its unique $G_2$ structure with parallel torsion.
\end{conjecture}
Equivalently, we expect the classes of closed $\phi$ with extremally
pinched Ricci curvature and those with parallel intrinsic torsion to
coincide on compact manifolds.

In the last section we present some examples intended to elucidate the
difficulty of making conclusions about torsion from vanishing of
curvature components, for general $G_2$ structures.  As a starting
point warped products and cohomogeneity-one $G_2$ manifolds with
exactly one non-vanishing curvature component $(W_{77})$ and
non-trivial torsion are given.  These realize every possible torsion type
allowed by the vanishing of the scalar curvature, with one exception.
\begin{ack} 
  The authors wish to thank Robert Bryant, Andrew Swann, Simon
  Salamon, Ilka Agricola and Thomas Friedrich for inspiration and for
  making their help available.  R.~C. thanks Ilka Agricola and Thomas
  Friedrich for their support, patience and the stimulating work
  environment they provide, and the Volkswagen Foundation and Yat Sun
  Poon for support. S.~I. thanks the Max-Plank-Institut f\"ur
  Mathematik, Bonn for the support and the excellent research
  environment. S.~I. is a Senior Associate to the Abdus Salam ICTP,
  2003-2009.
\end{ack}

\section{The fundamental three-form}

Let $(V,\inp{\cdot}{\cdot})$ be an $n$-dimensional Euclidean vector
space.  An $n$-form $\vol$ on $V$ is said to be a volume form for the
given inner product if $\abs{\vol}=1$ with respect to the inner
product induced on the exterior algebra $\Lambda^* V^*$.  Write
$i_v\colon \Lambda^pV^*\to\Lambda^{p-1}V^*$ for the interior product.
Suppose now that $n=7$ and that $\phi$ is a three-form on $V$ such
that
\begin{equation}
  \label{eq:1}
  i_u\phi\bw i_v\phi\bw\phi=6\inp{u}{v}\vol
\end{equation}
for some positive definite inner product $\inp{\cdot}{\cdot}$ and
volume form $\vol$.  Then $\phi$ is non-degenerate in the sense that
$v\mapsto i_v\phi$ is injective.  It follows that the isotropy group
of $\phi$ is the simple Lie group $G_2$, see
e.g.~\cite{math.DG/0010054}, and a basis of orthonormal one-forms $e^i$
for $V^*$ may
be chosen such that
\begin{equation}
  \label{eq:2}
  \begin{split}
    \phi &= e^{127} + e^{347} + e^{567} + e^{135} - e^{245} - e^{146} -
    e^{236},\\
    \Hodge\phi &= e^{1234} + e^{3456} + e^{5612} - e^{2467} + e^{1367}
    + e^{2357} + e^{1457}.
  \end{split}
\end{equation}
Here $e^{ijk}$ is short-hand for $e^i\bw e^j\bw e^k$ and so on.  From
equation~\eqref{eq:1} it is clear that $G_2$ is a closed subgroup of
$\SO(7)$.  Via the inner product the Lie algebra $\lie g_2$ of $G_2$
may be identified with the $14$ dimensional subspace of $\Lambda^2V^*$
complementary to the span of $\{i_u\phi \colon u\in V\}$.  One says
that a three-form satisfying equation~\eqref{eq:1} for some positive
inner product and volume form on $V$ is a \emph{$G_2$ three-form} or
\emph{fundamental three-form of $G_2$}.  The inner product and volume
form so defined are said to be \emph{associated to the three-form}.
Alternatively, fixing an inner product and volume form, any three-form
satisfying the relation~\eqref{eq:1} is called \emph{compatible} with
the metric and given orientation. A basis of one-forms $\{e^i\}$ over
a vector space $V$ for which a $G_2$ three-form $\phi$ has the
expression~\eqref{eq:2} is called \emph{$G_2$-adapted}.

A $G_2$ three-form $\phi$ induces a splitting of the exterior algebra
$\Lambda^* V$.  We have equivariant projections
$p^r_d\colon\Lambda^rV^*\to\Lambda^rV^*$ given by
\begin{equation}
  \label{eq:4}
  \begin{array}{c}
    p^2_{7}(\alpha)=\tfrac13(\alpha+\Hodge(\alpha\bw\phi),\qquad
    p^2_{14}(\alpha)=\tfrac13(2\alpha-\Hodge(\alpha\bw\phi),\\
    p^3_{1}(\beta)=\tfrac17\Hodge(\Hodge\phi\bw\beta)\phi,\qquad
    p^3_{7}(\beta)=\tfrac14\Hodge(\Hodge(\phi\bw\beta)\bw\phi),\qquad
    p^3_{27}(\beta)=\beta-(p^3_{1}+p^3_{7})(\beta).
  \end{array}
\end{equation}
A subscript $d = 1,~7,~14$, or $27$ indicates the dimension of the
images, denoted by $\Lambda^r_d$, of the corresponding projections.
These are all irreducible representations of $G_2$.  Projections for
$r>3$ are obtained by composing with the Hodge star operator
$\Hodge\colon\Lambda^r\to\Lambda^{7-r}$, $p^r_d:=\Hodge\circ
p^{7-r}\Hodge$.

On several occasions in what follows we shall come across
representations that do not occur as subspaces in the exterior algebra
of the standard representation.  We fix the notation for these as
follows.  Choose a system of positive roots for $\lie g_2$ such that
the standard representations has highest weight $(1,0)$ and the
adjoint representation $(0,1)$.  We write $V^{(\mu_1,\mu_2)}_d$ for
(the isomorphism class of) the irreducible representation with highest
weight $(\mu_1,\mu_2)$ and dimension $d$.  So the standard
representation is $V_7^{(1,0)}$, the adjoint representation $\lie g_2$
is $V^{(0,1)}_{14}$, while the space of traceless symmetric tensors is
$V^{(2,0)}_{27}$.  When the dimension is sufficient to identify the
representation the superscript will be dropped.

The $27$ dimensional subspace $\Lambda^3_{27}$ is isomorphic the space
of traceless symmetric tensors over $V$.  This isomorphism may be
given explicitly as the restriction of 
\begin{equation}
  \label{eq:5}
  \lambda_3(e\otimes e):=e\bw (e\hook\phi)
\end{equation}
to tracefree tensors.  An map from $\Lambda^3$ to $2$-tensors is given
by contracting an arbitrary three-form with the fundamental form over
two indices
\begin{equation*}
  \sigma(\alpha)(u,v)=\inp{i_u\phi}{i_v\alpha}.
\end{equation*}
The two-tensor $\sigma(\alpha)$ is a symmetric tensor only when
$p^3_7(\alpha)$ vanishes, and tracefree only if $p^3_1(\alpha)=0$.
Note that $\lambda_3(g)=3\phi$.  For a symmetric tensor $h$ with zero
trace one has the simple relation
\begin{equation}
  \label{eq:7}
  \lvert{\lambda_3(h)}\rvert^2=2\norm{h}^2.
\end{equation}
A few words on how identities such as~\eqref{eq:7} and~\eqref{eq:4}
are verified, given that these techniques are well-established.  All
identities here represent relations between maps to or from an
irreducible representation.  Then Schur's Lemma ensures that any two
such maps must be equal up to a constant multiple.  It is then
sufficient to evaluate left- and right-hand sides on a \emph{test
  element}.  Examples of such may be provided as follows.  Choose a
$G_2$ adapted basis $e^i$ of $V^*$.  Set
$\omega:=e^{12}+e^{34}+e^{56}$ and $\psi^+:= e^{135} - e^{245} -
e^{146} - e^{236}$.  Then $\phi=\omega\bw e^7 + \psi^+$ and
$\omega\in\Lambda^2_7$, $e^{12}-e^{34}\in\Lambda^2_{14}$ are test
elements in $\Lambda^2V^*$.  In degree~$3$ we have
$\phi\in\Lambda^3_1$ while $\psi^-:=- e^{246} + e^{136} + e^{235} +
e^{145}\in\Lambda^3_7$ and $4\omega\bw e^7-3\psi^+\in\Lambda^3_{27}$.
More test elements with an application of Schur's Lemma are given in
the proof of the following Lemma.
\begin{lemma}\label{lem:alg}
  Let $V$ denote $\mathbb R^7$ equipped with $G_2$ three-form
  $\phi$~\eqref{eq:2} and associated metric $g$.  Let $\bw_3\colon
  V^*\otimes \Lambda^2 V^*\to \Lambda^3 V^*$ be the map given by
  the wedge product $\bw_3(\alpha\otimes\beta) = \alpha\bw\beta$.
  The tensor product $V^*\otimes \Lambda^2_{14}$ decomposes as
  \begin{equation*}
    V^*\otimes \Lambda^2_{14} \cong
    V^{(1,1)}_{64}+V^{(2,0)}_{27}+V^{(1,0)}_7 
  \end{equation*}
  and the restriction $\bw_3\rvert\colon V^*\otimes \Lambda^2_{14}
  \to \Lambda^3 V^*$ has kernel $V^{(1,1)}_{64}$ and cokernel
  $\Lambda^3_1$.  Moreover, the identity
  \begin{equation*}
    7\norm{\gamma}^2=\norm{\bw_3(\gamma)}^2
  \end{equation*}
  holds for every $\gamma$ in the $27$ dimensional irreducible
  submodule of $V^*\otimes \Lambda^2_{14}$.
\end{lemma}
\begin{proof}
  First, since $V^*$ and $\Lambda^2_{14}$ are non-isomorphic
  representations the decomposition of their tensor product contains
  no trivial summand.  So it is clear that the cokernel of
  $\bw_3\colon V^*\otimes \Lambda^2_{14}\to \Lambda^3 V^*$ must
  contain $\Lambda^3_1$.  Furthermore, as $V^{(2,0)}_{27}$ is real and
  irreducible, there is, up to scale, precisely one invariant map
  $S^2(V^{(2,0)}_{27})\to\mathbb R$.  By Schur's Lemma, there exists a
  constant $c$ so that the relation $c\norm{\gamma}^2 =
  \norm{\bw_3(\gamma)}^2$ holds for all $\gamma\in V^{(2,0)}_{27}$.
  
  Let $\{e^i\}$ be a $G_2$ adapted basis.  Write $\pi\colon V^*\otimes
  \Lambda^2 V^*\to V^*\otimes \Lambda^2_{14}$ for the orthogonal
  projection $\alpha\otimes\beta \mapsto \alpha\otimes
  p_{14}^2(\beta)$.  Then $\pi(e^i\otimes e^{i7})$ provides a test
  element in a submodule of $V^*\otimes \Lambda^2_{14}$ isomorphic to
  $V^*\cong V^{(1,0)}_7$ and one may calculate $\bw_3(\pi(e^i\otimes
  e^{i7}))=-\psi^-\in\Lambda^3_7$.  This shows that
  $\coker(\bw_3\rvert)$ is no bigger than $\Lambda^3_{27} +
  \Lambda^3_1$.
  
  Set $\gamma' := e^7\otimes(e^{12} - e^{34})\in V^*\otimes
  \Lambda^2_{14}$.  Then $\bw_3(\gamma') =e^{127} - e^{347} \in
  \Lambda^3_{27}$ which proves that $\coker(\bw_3\rvert)=\Lambda^3_1$
  and also shows that $V^*\otimes \Lambda^2_{14}$ contains irreducible
  submodules isomorphic to $V_{27}$ and $V$.  The decomposition now
  follows by noting that the dimension of the Cartan product
  $V^{(1,1)}$ inside $V^*\otimes \Lambda^2_{14}$ is $64$.  It is then
  clear that $\ker(\bw_3\rvert)\cong V_{64}^{(1,1)}$.
  
  All we now need is to find a test element in $V^{(1,1)}_{64}\subset
  V^*\otimes \Lambda^2_{14}$.  Since $p^3_7(\bw_3(\gamma'))=0$,
  $\gamma'$ itself must lie the submodule isomorphic to
  $V_{64}^{(1,1)}+V_{27}^{(2,0)}$.  Composing the inclusion
  $i\colon\Lambda^3V^* \hookrightarrow V^*\otimes \Lambda^2 V^*$ with
  the projection $\pi$ we obtain $\gamma'' : = \pi(i(\bw_3(\gamma')))
  - \gamma'$.  It is easy to check that $\inp{\gamma''}{\gamma'} = 0,~
  \norm{\gamma''}^2 = \tfrac{16}3,~ \norm{\gamma'}^2 = 4$ in the
  tensor norm, and $\bw_3(\gamma'')=\tfrac43\bw_3(\gamma')$.  We then
  have $4\gamma''-3\gamma'\in V_{64}^{(1,1)}$ and $\gamma := \gamma''
  + \gamma'\in V_{27}^{(2,0)}$.  Evaluating the norms of $\gamma$ and
  $\bw_3(\gamma)$ completes the proof.
\end{proof}

When modules are not irreducible but maps between them still $G_2$,
calculations on the components of the tensors are made easier by the
fact that $\phi$ generates the space of invariant tensors.  This
guarantees that relations between contractions of metric $g,~\phi$ and
$\Hodge\phi$ - spelled out in~\cite{Br1} and also exploited
in~\cite{CI} - do exist.  For ease of reference we recall these
identities here.  Let $\phi$ be a $G_2$-three-form and $\Hodge\phi$
its dual four-form via the associated metric and orientation.  Write
$\phi_{ijk}$ for the components of $\phi$ and $\phi_{ijkl}$ for the
components of $\Hodge\phi$ with respect to a basis $e^i$ of one-forms
on $V^*$.  Then one has the identities, see~\cite{Br1}
\begin{align}
  \label{eq:8}
  \phi_{ipq}\phi_{pqj} &= 6\delta_{ij},\\
  \label{eq:9}
  \phi_{ijp}\phi_{pkl} &=
    \delta_{ik}\delta_{jl}-\delta_{jk}\delta_{il} + \phi_{ijkl},\\
  \label{eq:10}  
  \phi_{ijpq}\phi_{pqkl} &=
    4(\delta_{ik}\delta_{jl}-\delta_{jk}\delta_{il}) + 2\phi_{ijkl},\\
  \label{eq:11}  
  \phi_{ipq}\phi_{pqjk} &= 4\phi_{ijk},\\
  \label{eq:12}  
  \phi_{ijp}\phi_{pklm} &= \delta_{ik}\phi_{jlm} -
  \delta_{jk}\phi_{ilm} + \delta_{il}\phi_{jmk} -
  \delta_{jl}\phi_{imk} + \delta_{im}\phi_{jkl} -
  \delta_{jm}\phi_{ikl},
\end{align}
where repeated indices here and below indicate that a summation is
taking place.  The identities listed here are valid only when the
basis chosen is orthonormal.  For a general basis one must replace the
$\delta_{pq}$'s with the components $g_{pq}$ of the associated metric
in this basis and simple summations must be replaced with contractions
with the inverse metric.  This means, for instance, that the first
identity becomes $\phi_{ipr}g^{pq}g^{rs}\phi_{qsj} = 6g_{ij}$. 

\section{Torsion of a $G_2$ structure}
 
A \emph{$G_2$ structure or $G_2$ three-form} on a $7$ dimensional
manifold $M$ is a three-form $\phi$ such that for any two vector
fields $X,Y$
\begin{equation}\label{eq:15}
  i_X\phi\bw i_Y\phi\bw\phi=6g(X,Y)\vol(g),
\end{equation}
where $g$ is a Riemannian metric and $\vol(g)$ is a volume element for
$g$.  Fixing the three-form, we say that $g$ and $\vol(g)$
satisfying~\eqref{eq:15} are the metric and orientation
\emph{associated} to $\phi$.  A different viewpoint is offered by
fixing a metric $g$ and a orientation.  Then a three-form $\phi$ is
called \emph{compatible} with this choice when~\eqref{eq:15} holds.

A $G_2$ structure defines reductions of the bundle $LM$ of linear
frames $\sigma\colon\mathbb R^7\to T_xM$ of $M$ first to the bundle of
orthonormal frames $F(M,g)$ such that $\sigma^* g$ is the standard
metric on $V=\mathbb R^7$ and then to the bundle $P$ of $G_2$ adapted
frames such that $\sigma^*\phi$ is the standard $G_2$ structure on $V$
given by \eqref{eq:2}.  Any such reduction of the bundle of
orthonormal frames with structure group $G$ carries a unique $G$
connection $\Nt$ called \emph{the canonical connection of $P$}.  This
has the property that \emph{the intrinsic torsion $\xi$} defined by
$\xi:=\LC-\Nt$ of $P$ is a one-form with values in $\lie
g^{\perp}\subset \lie{so}(7)$.  Here and everywhere else in this paper
$\LC$ denotes the Levi-Civita connection of a given metric $g$.
Identifying tensor bundles $T^{(p,q)}M$ with the corresponding
associated bundles $P\times_G V^{(p,q)}$ allows to identify tensor
fields with equivariant functions $P\to V^{(p,q)}$.  When $W$ is a
$G$-representation, we adopt the convenient but somewhat improper
notation $\gamma\in W$, meaning $\gamma$ is a section of the bundle
associated to $W$.  Given a $G$-invariant tensor $\gamma$ on $V$ it
makes sense to define a section $\gamma \in W$ by setting $p\mapsto
\gamma$ for all $p\in P$.  This tensor will be parallel with respect
to $\Nt$.  In this way, equivariant linear maps between
representations $c\colon W\to W'$ of $G$ give rise to maps of
associated bundles commuting with the covariant differentiation by
$\Nt$
\begin{equation}
  \label{eq:13}
  \Nt_X(c(\eta))=c(\Nt_X\eta), \qquad X\in V,~\eta\in W.
\end{equation}
Let $.$ denote the induced action of a Lie algebra $\lie g$ on
representations $V$ of $G$.  Extending this to associated bundles we
apply the short-hand notation $\LC\gamma = \Nt\gamma + \xi.\gamma$.
With the convention that $\bw_p\colon
V^*\otimes\Lambda^pV^*\to\Lambda^{p+1}V^*$ is the map
$\alpha\otimes\beta\to\alpha\bw\beta$ we define
\begin{equation*}
  d^{\Nt}\beta:= \bw_p(\Nt\beta) = d\beta-\bw_p(\xi.\beta).
\end{equation*}

For further details on $G$ structures see for instance \cite{Salamon}.

The complement $\lie g_2^{\perp}$ of $\lie g_2$ in $\so(7)$ is
isomorphic to the standard representation $V=V_7$ of $\Lie G_2$.  So
the decomposition of the intrinsic torsion $\xi$ of a $G_2$ structures
follows from the splitting
\begin{equation}
  \label{eq:14}
  V\otimes V = S^2V+\Lambda^2V = S^2_0V + \mathbb R + \Lambda^2_{14} +
  \Lambda^2_{7}. 
\end{equation}
We write $\xi_d$ for the projection of $\xi$ to the $d$-dimensional
subspace of $V^*\otimes\lie g_2^\perp$ corresponding to the
decomposition~\eqref{eq:14} $\xi=\xi_{27}+\xi_1+\xi_{14}+\xi_7$.

\subsection{Derivatives of the fundamental three-form}
\label{sec:deriv-fund-three}
The \emph{torsion components} of a $G_2$ structure are differential
forms $\tau_p\in\Omega^p(M)$ such that~\cite{FG}
\begin{equation}
  \label{eq:16}
  \begin{array}{c}
    d\phi = \tau_0\Hodge\phi + 3\tau_1\bw\phi+\Hodge\tau_3,\\
    d\Hodge\phi = 4\tau_1\bw\Hodge\phi + \tau_2\bw\phi.
  \end{array}
\end{equation}
This pair of equations are the \emph{structure equations} for the
$G_2$ form $\phi$.

The torsion type or Fern\'andez-Gray class of a $G_2$ structure is
determined by the vanishing of torsion components.  The original
enumeration of the torsion types of Fern\'andez and Gray gives the
correspondences $\tau_0\leftrightarrow 1$, $\tau_1\leftrightarrow 4$,
$\tau_2\leftrightarrow 2$, and $\tau_3\leftrightarrow 3$.  For
instance, a three-form has (Fern\'andez-Gray) type $1+3$ if
$\tau_2=0=\tau_1$ and strict, or proper, type $1+3$ if
$\tau_2=0=\tau_1$, but $\tau_0\not\equiv0$ and $\tau_3\not\equiv 0$.
If all components are zero then $\phi$ is parallel.

The torsion $\tau$ and intrinsic torsion may be related explicitly as
follows.  Let $\bar\xi$ denote the two-tensor obtained by
identification of $\lie
g_2^{\perp}\simeq\Lambda^2_{7}$ with $\Lambda^1_7\simeq V^*$.  In an orthonormal
frame $e_i$ this is the contraction
\begin{equation}
  \label{eq:17}
  \bar\xi=\xi_{ipq}\phi_{pqj}e^i\otimes e^j,
\end{equation}
where $\xi_{ijk}=g(\xi_{e_i}e_j,e_k)$.  The components of the
intrinsic torsion may be recovered from $\bar\xi$ by the relation
\begin{equation}
  \label{eq:3}
  \xi_{ijk}=\tfrac16\bar\xi_{ip}\phi_{pjk}.  
\end{equation}
In fact, still working in an orthonormal frame $e_i$ one then has a
nice expression for the covariant derivative of $\phi$:
$(\LC_i\phi)_{jkl}=-\tfrac12\bar\xi_{ip}\phi_{pjkl}$, which leads to
the relations
\begin{equation}
  \label{eq:19}
  \begin{array}{ll}
    \bar\xi_1=-\tfrac12\tau_0 g,&
    \bar\xi_7=2\Hodge(\tau_1\bw\Hodge\phi),\\
    \bar\xi_{14}=\tau_{2},&
    \bar\xi_{27}=\sigma(\tau_3).
  \end{array}
\end{equation}
These may also be found in~\cite{Kar}.  An example of an application
of the identity~\eqref{eq:13} is the component-wise relation of
covariant derivatives of $\bar\xi$ and $\xi$, $\Nt_i\xi_{jkl} =
\Nt_i\bar\xi_{jp}\phi_{pkl}$.

\section{Curvature of $G_2$ manifolds}
\label{sec:curv-g_2-manif}
For a $G$ structure on $(M,g)$ the Riemannian curvature tensor may be
given the following expression in terms of the canonical connection
$\Nt$ and the intrinsic torsion $\xi$:
\begin{equation}\label{eq:20}
  \RC = \Rt+(\Nt\xi)+(\xi^2).
\end{equation}
Here, $\Rt\in\Lambda^2\otimes\lie g$ is the curvature of the canonical
connection, $(\Nt\xi)\in\Lambda^2V^*\otimes\lie g^{\perp}$ is defined
by $(\Nt\xi)_{X,Y} := (\Nt_X\xi)_Y - (\Nt_Y\xi)_X$ and $(\xi^2) \in
\Lambda^2V^*\otimes\so(7)$ is the tensor $(\xi^2)_{X,Y}Z =
\xi_{\xi_YX-\xi_XY}Z + [\xi_X,\xi_Y]Z$. Let $\mK(\lie g)$ be the space
of algebraic curvature tensors with values in $\lie g\subset \so(n)$
and fix $\mK:=\mK(\so(n))$.  Write $\mK=\mK(\lie g)\oplus\mK(\lie
g)^{\perp}$ with respect to $g$.  The first statement of the following
proposition is then an easy consequence of the first Bianchi identity
for $\RC$.  The second statement is a consequence of the first and a
theorem of Ambrose and Singer~\cite{MR0102842}, see~\cite{MR2114426}.
\begin{proposition}\label{thm:curv}
  Let $P$ be a $G$ structure on a Riemannian manifold $M$ with metric
  $g$.  Then the components of the Riemannian curvature in $\mK(\lie
  g)^{\perp}$ are determined by the components of the covariant
  derivative of the intrinsic torsion $\Nt\xi$ and the tensor
  $(\xi^2)$.  
  
  If $g$ is complete, $\Nt\xi=0$ and $\mK(\lie g)=\{0\}$ then $(M,g)$
  is locally isometric to a homogeneous space and the universal cover
  of $(M,g)$ is globally homogenous.
\end{proposition}
\begin{remark}
  Note that one may equally well express the Riemannian curvature as
  $\RC = \Rt+(\LC\xi)-[\xi^2]$, where $[\xi^2]_{X,Y}:=[\xi_X,\xi_Y]$.
  In the almost Hermitian case $\lie g=\un(m)$ and $\lie
  g^{\perp}=\real{\Lambda^{(2,0)}}$.  It then follows that $[\xi^2]
  \in \Lambda^2V^* \otimes ((\real{\Lambda^{(2,0)}} \otimes
  \real{\Lambda^{(2,0)}}) \cap
  \so(2m))\subset\Lambda^2V^*\otimes\un(m)$.  The bracket
  $\real{\cdot}$ sends a complex representation to its underlying real
  representation.  Based on this it was argued in~\cite{SalamonFF}
  that the components of Riemannian curvature in $\mK^\perp$ are
  determined by the components of $\LC\xi$.  However, this does not
  carry over to the $G_2$ setting simply because $[\lie g_2^\perp,\lie
  g_2^\perp]\not\subset\lie g_2$ as one easily verifies.
\end{remark}
For $G_2$ structures the decomposition of $\mK$ is easily obtained.
First, using standard techniques of representation theory we have
\begin{equation}
  \label{eq:21}
  \begin{array}{ccc}
    S^2(\lie g_2^\perp)=
    V_{27} + V_1,&S^2(\lie g_2) =
    V^{(2,0)}_{77} + V_{27} + V_1, &\lie g_2\odot\lie g_2^{\perp}
    =
    V_{64} + V_{27} + V_7,\\ 
    S^2(\lie g_2^\perp) \cap \mK = 0, & S^2(\lie g_2)\cap \mK =
    V^{(2,0)}_{77}, & (\lie g_2\odot\lie g_2^{\perp}) \cap \mK =
    V_{64}. 
  \end{array}
\end{equation}
where the convention $V\odot W := (V\otimes W+W\otimes V) \cap
S^2(V+W)$ for vector spaces $V$ and $W$ is used.  The decomposition of
$S^2(\so(7))$ follows from \eqref{eq:21}.  Let $\bian\colon\Lambda^2
V^*\otimes
\End(V)\to\Lambda^3V^*\otimes V$ be defined as $(\bian r)(X,Y,Z) :=
r(X,Y)Z + r(Y,Z)X + r(Z,X)Y$ so that $\mK=\ker(\bian)$.  Using that
$\bian$ restricted to $S^2(\so(7))$ maps surjectively onto $\Lambda^4$
one obtains $\mK = V^{(2,0)}_{77} + V_{64} + 2V_{27} + V_1$, as
observed in~\cite{Br1}.  Comparing to the $\so(7)$ decomposition of
$\mK$, see e.g.~\cite{Besse:Einstein}, $\mK=\mW + \mR_0 +\mS$ where
$\mR_0\cong S^2_0V^*,~\mS\cong\mathbb R$, we see that the space of
algebraic Weyl tensors $\mW$ on a $G_2$ manifold decomposes as
\begin{equation*}
  \mW = \mW_{77}+\mW_{64}+\mW_{27},
\end{equation*}
where $\mW_{77} := \mK\cap S^2(\lie g_2) \cong V^{(2,0)}_{77}$,~
$\mW_{64} := \mK \cap (\lie g_2\odot\lie g_2^{\perp}) \cong
V^{(1,1)}_{64}$ and $\mW_{27} := \mW \setminus (\mW_{77} + \mW_{64})
\cong \Lambda^3_{27}\cong S^2_0 V_7$.  This, at least from the point
of view of splitting the space of algebraic curvature tensors $\mK$ in
to $G_2$ irreducible subspaces, gives the \emph{decomposition of the
  Riemannian curvature tensor of a $G_2$ manifold.}  However, it is
sometimes useful to have a more explicit description of these
submodules.  To attain this, we first need to do a little more linear
algebra.

So let for the moment $\phi,g$ be the standard $G_2$ structure on
$V_7=V=\mathbb R^7$.  Let $r_g$ be the usual Kulkarni-Nomizu product
viewed as an $\SO(7)$ equivariant map $S^2V^*\to S^2(\Lambda^2V^*)$,
\begin{multline}\label{eq:22}
  r_g(h)(x,y,z,w) := (h\ovee
  g)(x,y,z,w)=\\h(y,z)g(x,w)-h(x,z)g(y,w)
  +h(x,w)g(y,z)-h(y,w)g(x,z).
\end{multline}
This takes values in $\mK$, as one may easily verify.  Then a $G_2$
equivariant map $r_\phi$ also from $S^2V^*$ to $\mK$ is
\begin{equation*}
  r_\phi(a_1\odot a_2):=(a_1\hook\phi)\odot(a_2\hook\phi) -
  \tfrac13\bian\left((a_1\hook\phi)\odot(a_2\hook\phi)\right),
\end{equation*}
Here and elsewhere, $a\odot b :=a\otimes b+b\otimes a$.  The Bianchi
map must of course be composed with the proper musical morphisms.
Contractions going in the opposite direction may be given as
\begin{gather*}
  c^g(r)(u,v):= r(u,e_i,e_i,v),\\
  \intertext{where $e_i$ is an orthonormal basis. This is just the
    usual Ricci contraction.  Using the isomorphism
    $S^2(\Lambda^2V)\cong_g S^2(\Lambda^2V^*)$ we set}
  c^\phi(r)(u,v):=4r(u\hook\phi,v\hook\phi).
\end{gather*}
The first equations of~\eqref{eq:23} below,
correspond to the result $c^g(h\owedge g)=(n-2)h+\tr_g(h)g$ of taking
the Ricci contraction of a Kulkarni-Nomizu product,
see~\cite{Besse:Einstein}.  These and the remaining equations may be
verified by a calculation in an orthonormal basis.
\begin{gather}
  \label{eq:23}
  \begin{array}{lll}
    (c^g\circ r_g)\lvert_{S^2_0V^*} = 5, & (c^g\circ
    r_\phi)\lvert_{S^2_0V^*} = 1, & (c^g\circ r_g)(g)=12g,\\
    (c^\phi\circ r_g)\lvert_{S^2_0V^*} = 4,
    & (c^\phi\circ
    r_\phi)\lvert_{S^2_0V^*} = \tfrac{92}3, & (c^\phi\circ r_g)(g)=-24g.
  \end{array}
\end{gather}
In analogy with the characterization $\mR_0=\{r_g(h)\colon h\in
S^2_0V^*\}$, the space $\mW_{27}$ may be described as
\begin{equation*}
  \mW_{27}:=\left\{r_g(h)-5r_\phi(h)\colon h\in S^2_0V^*\right\}.
\end{equation*}
At last consider the projections $ P^{\lie g_2}\colon S^2(\so(7))\to
S^2(\lie g_2),~ P^\odot\colon S^2(\so(7))\to \lie g_2\odot\lie
g_2^{\perp}$.  These maps may be given closed form expression in terms
of the projections $p^2_d,~d=7,~14$ of~\eqref{eq:4}.

\subsection{Ricci curvature of $G_2$ manifolds}
\label{sec:ricci-curvature-g_2}
Let $(M,\phi)$ be a $G_2$ manifold with associated metric $g$ as
usual.  Define the \emph{$\phi$-Ricci tensor}\footnote{This is
  called the $*$-Ricci tensor in~\cite{CI}} as 
\begin{equation}
  \label{eq:25}
  \Ric^\phi(X,Y):=c^\phi(\RC),
\end{equation}
and write $\Ric^g=c^g(\RC)$ and $s_g=\tr_g(\Ric^g)$ for the Ricci and
scalar curvatures of $g$.  The identities~\eqref{eq:23} show that
$\tr_g(\Ric^\phi)=-2s_g$, see also~\cite{CI}.  As usual, a subscript
$0$ indicates the traceless part of a symmetric tensor.  

The isomorphism \eqref{eq:5} has striking consequences. Most
importantly, the covariant derivatives of $\xi_7,\xi_{14}$ and
$\xi_{27}$ all have precisely one component in a $27$ dimensional
irreducible subspace of $V^*\otimes V^*\otimes\lie g_2^{\perp}$.  Each
of these may be identified with $27$ dimensional components of
suitable exterior derivatives.  Similarly, each $27$ dimensional
component of the `algebraic' parts $\xi_d\odot\xi_{d'}$ of the
Riemannian curvature has an equivalent expression in the exterior
algebra.  This was used in~\cite{Br1} to obtain an expression for the
Ricci curvature of a $G_2$ manifold.  When $\tau_2=0$, Ricci is given
in terms of covariant derivatives of the skew-symmetric torsion
in~\cite{SITF} and a formula for the Ricci tensor in terms of the
covariant derivatives of the intrinsic torsion was very recently
presented in~\cite{Kar}.  All this is subsumed by introducing a more
general symmetric two-tensor, built $G_2$ equivariantly from the
Riemannian curvature, as follows.  For $k=(k_1,k_2)\in \mathbb R^2$
write
\begin{equation}
  \label{eq:27}
  \Ric^k_0:=k_1\Ric^g_0+k_2\Ric^\phi_0.
\end{equation}
We shall call this tensor the \emph{generalized Ricci tensor} and say
that a $G_2$ manifold $(M,\phi)$ is \emph{generalized Einstein} if
$\Ric^k_0=0$ for some $k=(k_1,k_2)$.  
\begin{remark}
  Since rescaling $\Ric^k_0$ does not affect the generalized Einstein
  equation there is of course only an $\mathbb RP(1)$-worth of such
  constraints, really.
\end{remark}

As $S^2_0V$ has multiplicity two in $\mK$, any traceless symmetric two
tensor on a $G_2$ manifold must equal $\Ric^k_0$ for some $k\in\mathbb
R^2$.  Before giving the formula we need to demonstrate how the
components $\tau_2\odot\tau_3$ and $\tau_3\otimes\tau_3$ determine
three-forms.  For this, choose a local $G_2$ adapted frame
$(e_1,\dots,e_7)$.  For a two-form $\alpha$ and a three-form $\beta$
set
\begin{gather*}
  [\alpha\odot\beta]:=\sum_k i_{e_k}\alpha\bw i_{e_k}\beta,\quad
  [\beta^2]^A:=\sum_k \Hodge(i_{e_k}\beta\bw i_{e_k}\beta),\quad
  [\beta^2]^B:=\sum_k ((i_{e_k}\phi)\hook\beta)\bw i_{e_k}\beta.
\end{gather*}
These are independent of the chosen frame and so extend to smooth
contractions on $M$. The first two contractions are in fact $\SO(7)$
equivariant.  The fact that $V_7$ does not occur in the decomposition
of $S^2(V_{27})$ guarantees that for $\beta\in\Lambda^3_{27}$,
$[\beta^2]^A$ and $[\beta^2]^B$ belong to $\Lambda^3_1 +
\Lambda^3_{27}$.  Since there is a summand isomorphic to $V_7$ in
$\Lambda^2_{14}\otimes\Lambda^3_{27}$ it is not quite obvious that
$[\alpha\odot\beta]\in \Lambda^3_1+\Lambda^3_{27}$ should hold for
$\alpha\in\Lambda^2_{14}$ and $\beta\in\Lambda^3_{27}$, but this is
nevertheless true.

\begin{lemma}\label{prop:1}
  Let $(M,\phi)$ be a $G_2$ manifold.  Then 
  \begin{multline}
    \label{eq:28}
    \lambda_3\bigl(\Ric^k_0\bigr) =
    \bigl(-(5k_1+4k_2)d\Hodge(\tau_1\bw\Hodge\phi) +
    2(5k_1+4k_2)\tau_1\bw\Hodge(\tau_1\bw\Hodge\phi)\\
    - (k_1-4k_2)d\tau_2 + \tfrac12(k_1+2k_2)\Hodge(\tau_2\bw\tau_2) +
    (k_1+4k_2)\Hodge d\tau_3 + k_2[\tau_3^2]^A +
    \tfrac12k_1[\tau_3^2]^B\\
    -\tfrac12(k_1 - 4k_2)\tau_0\tau_3 + (k_1 - 4k_2)\tau_1\bw\tau_2 +
    (3k_1 - 4k_2)\Hodge(\tau_1\bw\tau_3) + 2k_2
    [\tau_2\odot\tau_3]\bigr)_{27}.
  \end{multline}
\end{lemma}
\begin{proof}
  Proposition~\ref{thm:curv} and the remarks above show that any symmetric
  two-tensor coming from contraction of $\RC$ with $\phi$ and $g$ must
  be a linear combination of the terms on the right-hand side where
  coefficients are determined entirely in terms of the linear algebra
  of $G_2$ and so in particular are independent of the underlying
  manifold.  Obtaining the given expression is then a matter of
  evaluating the left- and right-hand side on examples.
\end{proof}
For $k=(1,0)$ this is due to Bryant~\cite{Br1} and so is the scalar
curvature of a $G_2$ manifold: (see also~\cite{MR2006222,Kar}
\begin{equation}
  \label{eq:29}
  s^g=\tfrac{21}8\tau_0^2 + 12\delta\tau_1 + 30\abs{\tau_1}^2
  -\tfrac12\abs{\tau_2}^2-\tfrac12\abs{\tau_3}^2. 
\end{equation}

We also introduce the so-called $\mW$-Ricci-tensor:
\begin{equation*}
  \Ric^\mW:=\tfrac1{20}\left(4\Ric_0^g-5\Ric^\phi_0\right).
\end{equation*}

Lemma~\ref{prop:1} gives 
\begin{corollary}  \label{co:conf-inv}
  The symmetric traceless tensor $\Ric^{\mW}$ is a conformal
  invariant of a $G_2$ structure.
\end{corollary}
\begin{proof}
  Under a conformal change $\phi\to \tilde\phi = e^{3f}\phi$ the torsion
  is known to transform as $\tau = (\tau_0,\dots,\tau_3) \to
  \tilde\tau = (e^{-f}\tau_0,\tau_1+df,e^f\tau_2,e^{2f}\tau_3)$.
  Putting $\tilde\tau$ in to formula~\eqref{eq:28} a simple calculation
  verifies that with $k_1=4,~k_2=-5$ the generalized Ricci tensor
  merely rescales.
\end{proof}
\begin{remark}\label{rem:2}
  The subscript ${27}$ at the end of equation~\eqref{eq:28} refers to
  the image under the projection $p^3_{27}$.  It is rather pleasing to
  note that all possible contributions are present in the formula.
  Since there is a two parameter family of generalized Ricci tensors
  and a two parameter family of contributions from
  $\tau_3\otimes\tau_3$ it is not surprising that contributions from
  the chosen representatives $A$ and $B$ vanish for certain values
  $k$.  That those are precisely $(1,0)$ and $(0,1)$ appears to be
  pure coincidence.
  
  As noted in~\cite{Br1}, a generic $G_2$ structure has a
  two-parameter family of `canonical' $G_2$ connections
  $\nabla^{(s,t)}$.  The generalized Ricci-curvature defined above
  should therefore have an interpretation as the symmetric part of the
  contraction $c_g(R^{\nabla^{(s,t)}})$ of the curvature
  $R^{\nabla^{(s,t)}}$ of $\nabla^{(s,t)}$.  One clearly obtains
  different formulas for the Ricci-curvatures by expressing the
  exterior derivatives of the torsion components in terms of
  $d^{\nabla^{(s,t)}}$ instead.  We do this below for $\Nt$.
  
  Obtaining formulas for curvature components in the manner described
  above is also possible for other geometries.  In the case of
  $\SU(3)$-structures this has been done for the Ricci curvature
  in~\cite{math.DG/0606786}.
\end{remark}
A rather long but straightforward computation on test elements shows
that
\begin{multline*}
  d\Hodge(\tau_1\bw\Hodge\phi) = d^{\Nt}\!\Hodge(\tau_1\bw\Hodge\phi) +
  \tfrac12\Hodge(\tau_0\tau_1\bw\phi) +
  \tfrac83\tau_1\bw\Hodge(\tau_1\bw\Hodge\phi) - 2\abs{\tau_1}^2\phi\\
   + \tfrac13\tau_1\bw\tau_2 + \tfrac43(\tau_1\bw\tau_2)_7 +
  \tfrac23\Hodge(\tau_1\bw\tau_3) -
  \tfrac43\Hodge(\tau_1\bw\tau_3)_7,
\end{multline*}
\begin{multline}\label{eq:30}
  d\tau_2 = d^{\Nt}\!\tau_2 + \tfrac23\tau_1\bw\tau_2 +
  \tfrac83(\tau_1\bw\tau_2)_7 + \tfrac16\Hodge(\tau_2\bw\tau_2) +
  \tfrac16\abs{\tau_2}^2\phi\\
  - \tfrac16 [\tau_2\odot\tau_3] +
  \tfrac16\Hodge((\tau_2\hook\tau_3)\bw\phi),
\end{multline}
\begin{multline*}
  d\tau_3 = d^{\Nt}\!\tau_3 -\tfrac16\Hodge(\tau_0\tau_3) +
  \tau_1\bw\tau_3 - \tfrac83\Hodge(\Hodge(\tau_1\bw\tau_3)_7) -
  \tfrac16(\tau_2\hook\tau_3)\bw\phi\\
  - \tfrac16\Hodge[\tau_3^2]^A -\tfrac16\Hodge[\tau_3^2]^B +
  \tfrac16\abs{\tau_3}^2\Hodge\phi.
\end{multline*}
It is somewhat remarkable that there is no summand
$[\tau_2\odot\tau_3]$ in the last expression.  More surprises are in
store when these expression are inserted in formula~\eqref{eq:28}.
First, for $\beta\in\Lambda^3 V^*$ define
\begin{equation*}
  [\beta^2]^C:=[\beta^2]^A - 2 [\beta^2]^B.
\end{equation*}
Then Lemma~\ref{prop:1} is reformulated as
\begin{lemma}\label{prop:new1} 
Let $(M,\phi)$ be a $G_2$ manifold. Then
\begin{multline}\label{eq:new1}
  \lambda_3\bigl(\Ric^k_0\bigr) =
  \bigl(-(5k_1+4k_2)d^{\Nt}\!\Hodge(\tau_1\bw\Hodge\phi)
  -\tfrac23(5k_1+4k_2)\tau_1\bw\Hodge(\tau_1\bw\Hodge\phi)\\ -
  (k_1-4k_2)d^{\Nt}\!\tau_2 +\tfrac13(k_1+5k_2)\Hodge(\tau_2\bw\tau_2)
  + (k_1+4k_2)\Hodge d^{\Nt}\!\tau_3 - \tfrac16(k_1 -
  2k_2)[\tau_3^2]^C \\ -\tfrac23(k_1 - 2k_2)\tau_0\tau_3 -
  \tfrac43(k_1 + 2k_2)\tau_1\bw\tau_2 + \tfrac23(k_1 -
  4k_2)\Hodge(\tau_1\bw\tau_3) + \tfrac16(k_1 + 8k_2)
  [\tau_2\odot\tau_3]\bigr)_{27}.
\end{multline}
\end{lemma}
\begin{corollary}
  The components of the covariant derivative $\Nt\xi$ and the
  symmetric products $\xi_d\odot\xi_{d'}$ contribute to the traceless
  symmetric tensors $\Ric^g,\Ric^\phi$ and $\Ric^{\mW}$ according to
  Table~\ref{tab:g2-2}.
\end{corollary}
\begin{remark}
  Note that only one combination of a two parameter
  family of possible contributions from $\tau_3\otimes\tau_3$ is
  realized.
  
  A classification of $G_2$ structures with $\tau_1=\tau_2=0$ for
  which the torsion is parallel with respect to the unique $G_2$
  \emph{characteristic} connection determined by having three-form
  torsion: $T(X,Y)Z = - T(X,Z)Y$ was pursued
  in~\cite{math.DG/0604441}.  One corollary of this is that structures
  with torsion $\tau=\tau_3$, parallel with respect to the
  characteristic connection are never Einstein, at least when the
  stabilizer of the torsion is not Abelian.
  
  Note that when $d^{\Nt}\tau_3=0$ the Einstein equation for a $G_2$
  structure of type $1+3$ reduces to a quadratic equation
  $[\tau_3^2]^C+4\tau_0\tau_3=0$ for a torsion tensor
  $\tau=(\tau_0,\tau_3)$.  Moreover, if such a structure is Einstein
  then it is generalized Einstein for any value of $k=(k_1,k_2)$.  It
  is very likely that this system has non-trivial solutions.
\end{remark}
\begin{table}[tp]
\begin{tabular}{@{}l@{}ccc@{}}
\toprule
Tensor & \(\Ric_0^g\) & \(\Ric^{\mathcal W}\) & \(\Ric^{\phi}\) \\  
\midrule
\(\Nt\xi_7\)& \T & & \T \\
\(\Nt\xi_{14}\)& \T & \T  & \T \\
\(\Nt\xi_{27}\)&\T & \T & \T \\
\midrule
\(\xi_7\otimes\xi_7\)& \T & & \T \\
\(\xi_{14}\otimes\xi_{14}\)& \T & \T & \T \\
\((\xi_{27}\otimes\xi_{27})^C\)& \T & \T & \T\\
\midrule
\(\xi_{1}\odot\xi_{27}\)& \T &\T  & \T\\
\(\xi_{7}\odot\xi_{14}\)& \T& \T & \T \\
\(\xi_7\odot\xi_{27}\)& \T & \T & \T\\
\addlinespace[.75mm]
\(\xi_{14}\odot\xi_{27}\)& & \T & \T\\
\bottomrule
\end{tabular}
\caption{\(G_2\)-irreducible components of tensors contributing to
  Ricci curvatures.}
\label{tab:g2-2}
\end{table}
Now we can state one of our main results.
\begin{theorem}\label{thm:curv-decomp}
  Let $(M,\phi)$ be a $G_2$ manifold with associated metric $g$. Then
  the space of algebraic curvature tensors has the orthogonal
  splitting:
  \begin{gather*}
    \mK=\mW_{77}+\mW_{64}+\mW_{27}+\mR_0+\mS,
  \end{gather*}
  where the space of algebraic Weyl curvatures is
  $\mW=\mW_{77}+\mW_{64}+\mW_{27}$.  In terms of the scalar curvature,
  Ricci curvature and $\phi$-Ricci curvature the orthogonal
  projections to these subspaces are
  \begin{equation}\label{eq:26}
    \begin{split}
      S&=\tfrac1{84}s_g r_g(g),\\
      R_0&=\tfrac15r_g(\Ric^g_0),\\
      W_{27}&=\tfrac{3}{112}(r_g-5r_\phi)(\Ric^\mW),\\
      W_{64}&=P^\odot(W-W_{27}),\\ W_{77}&=P^{\lie g_2}(W-W_{27}).    
    \end{split}
  \end{equation}
  The three $G_2$ invariant components of the Weyl curvature are
  conformal invariants of the associated metric.
  
  In particular, $W_{27}=0$ exactly when $\Ric^\mW=0$ and the norm of
  the Riemannian curvature satisfies the equality
  \begin{equation*}
    \norm{\RC}^2 = \norm{W_{77}}^2 + \norm{W_{64}}^2 +
    \tfrac{15}{28}\lVert\Ric^\mW\rVert^2 + \tfrac45\norm{\Ric_0^g}^2 +
    \tfrac1{21}s_g^2. 
  \end{equation*}
\end{theorem}
\begin{proof}
  The first statement (which may also be found in~\cite{Br1}) is
  clear.  Write the Riemannian curvature as
  \begin{equation*}
    \RC = W_{77}+W_{64}+(r_g(h_1)-5r_\phi(h_1))+r_g(h_2)+Kr_g(g),
  \end{equation*}
  where $h_1$ and $h_2$ are traceless symmetric two-tensors, and apply
  the contractions $c^g$ and $c^\phi$ to obtain the first three
  equations in~\eqref{eq:26}.  The last two relations in~\eqref{eq:26}
  follow from the decompositions~\eqref{eq:21}.  The expression for
  the norm of the Riemannian curvature is now a matter of using the
  following relations
  \begin{equation*}
    \norm{r_g(h)}^2=20\norm{h}^2,\qquad
    \norm{r_\phi(h)}^2=\tfrac{92}3\norm{h}^2,\qquad
    \inp{r_\phi(h)}{r_g(h)}=4\norm{h}^2,
  \end{equation*}
  valid for any traceless symmetric two-tensor $h$, while
  $\norm{r_g(g)}^2=336$.   
  
  The conformal invariance of the three components $W_{77}$, $W_{64}$
  and $W_{27}$ follows from the conformal invariance of $W$ and
  invariance of orthogonal projection under rescaling of an inner product.
\end{proof}

\begin{table}[tp]
\begin{tabular}{@{}l@{}ccccc@{}}
\toprule
\lb{Tensor}&\multicolumn{4}{c}{Components in \(\mK\)} & \lb{Other
  components}\\\cmidrule(r){2-5} 
&\(V^{77}_{(0,2)}\) & \(V^{64}_{(1,1)}\) & \(V^{27}_{(2,0)}\) &
\(V^{1}_{(0,0)}\)&\\  
\midrule
\(\Nt\xi_1\)& & & & & \(V^7\)\\
\(\Nt\xi_7\)& & & \T & \T & \(V^7,~V^{14}\)\\
\(\Nt\xi_{14}\)& & \T & \T & & \(V^7\)\\
\(\Nt\xi_{27}\)& & \T & \T & & \(V^7,~V^{14},V^{77}_{(3,0)}\)\\
\midrule
\(\xi_1\otimes\xi_1\)& & & & \T & \\
\(\xi_7\otimes\xi_7\)& & & \T & \T & \\
\(\xi_{14}\otimes\xi_{14}\)& \T & & \T & \T & \\
\(\xi_{27}\otimes\xi_{27}\)& \T & \T& \(2\times\)\T & \T &
\(V^{182}_{(4,0)}\)\\
\midrule
\(\xi_1\odot\xi_7\)& & & & &\(V^7\)\\
\(\xi_1\odot\xi_{14}\)& & & & & \(V^{14}\)\\
\(\xi_{1}\odot\xi_{27}\)& & & \T & & \\
\(\xi_{7}\odot\xi_{14}\)& & \T & \T & & \(V^7\)\\
\(\xi_7\odot\xi_{27}\)& & \T & \T &  &\(V^7,~V^{14},V^{77}_{(3,0)}\)\\
\addlinespace[.75mm]
\(\xi_{14}\odot\xi_{27}\)& & \T & \T & &
\(V^7,~V^{14},~V^{189}_{(2,1)},V^{77}_{(3,0)}\)\\
\midrule
\(\Rt\)& \T & \T& \(2\times\)\T & \T & \(V^7,~V^{14},V^{77}_{(3,0)}\)\\
\bottomrule
\end{tabular}
\caption{\(G_2\)-irreducible components of tensors contributing to curvature.}
\label{tab:g2-1}
\end{table}
\begin{remark}\label{2order}
  Note that the components of $\Nt\xi$ and $\xi\odot\xi$ given in the
  right-most column of Table~\ref{tab:g2-1} correspond to the second
  order diffeomorphism invariants of $(M,\phi)$ not captured by the
  Riemannian curvature, see~\cite{Br1}.  Another interesting point
  made in~\cite{Br1} (and in~\cite{math.DG/0501062} for almost
  Hermitian structures) is that the first order identities for the
  torsion that derive from $d^2\phi=0=d^2\Hodge\phi$ are encoded in a
  subspace isomorphic to $2V_7+V_{14}$.  For $G_2$ these relations
  between invariants necessarily take their values in the complement
  of the space of algebraic curvature tensors inside the much larger
  space of diffeomorphism invariant polynomials in the derivatives of
  $\phi$ up to order two.  Andrew Swann has pointed out to us that
  these first order constraints on the torsion can be seen as coming
  from the fact that the cokernel of the restriction
  $\bian\colon\Lambda^2V^*\otimes\lie g_2\to\Lambda^3V^*\otimes V$ is
  isomorphic to precisely $2V_7+V_{14}$, and that the restriction
  $\bian\colon \Lambda^2V^*\otimes\lie g_2^{\perp} \to
  \Lambda^3V^*\otimes V$ is injective.
      
  For the examples in section~\ref{sec:1+4} and~\ref{sec:2} applying
  these principles show that all (polynomial) diffeomorphism
  invariants up to order $2$ are determined by the components of the
  Riemannian curvature tensor, a single one-form (which may be taken
  to be $d\tau_0$) and tensors, polynomial up to degree two in the
  torsion components.  This appears to be a distinguishing feature of
  the wider class of $G_2$-structures of type.
    
  We thank Robert Bryant for bringing the importance of the full space
  of second order diffeomorphism invariants of a $G_2$
  structure to our attention.
\end{remark}

\section{Curvature of $G_2$ structures of type $1+4$}
\label{sec:1+4}

Although the results of this section are well known, we take this
class $G_2$-structures as a first example to demonstrate how the
results of the previous sections can be used.  The inspection of the
structure equations~\eqref{eq:16} with $\tau_2=\tau_3=0$ yields:
$d\tau_1=0, ~d\tau_0 + \tau_0\tau_1 = 0$.  An easy argument now shows
that either $\tau_0\equiv 0$ and $\tau_1$ is closed, or $\tau_0$ is
never zero and $\tau_1$ exact: $\tau_1=-d\ln(\tau_0)$.  Therefore the
class $1+4$ of $G_2$-structures consists in $G_2$ structures $\phi$
locally conformally equivalent to a parallel structure and in those
globally conformally equivalent to a nearly parallel structure,
see~\cite{math.DG/0607487} for details.  From the entries in
Table~\ref{tab:g2-1} and Proposition~\ref{thm:curv} it is clear that
the Riemannian curvature of any such structure has $W_{64}=0$.  The
formulas of either of the Lemmas~\ref{prop:1} and~\ref{prop:new1} show
that $\Ric^\mW=0$ and so $W_{27}=0$ for type $1+4$ structures.  The
two remaining curvature components in $\mK(\lie g_2)^{\perp}$ are
determined by the torsion components as follows.
\begin{gather*}
  s_g= \tfrac{21}8\tau_0^2 + 12\delta\tau_1 + 30\abs{\tau_1}^2\qquad
  \lambda_3\left(\Ric_0^g\right)=\bigl(-5d^{\Nt}\!\Hodge(\tau_1\bw\Hodge\phi)
  -\tfrac{10}3\tau_1\bw\Hodge(\tau_1\bw\Hodge\phi)\bigr)_{27}.
\end{gather*}
\noindent {\bf Nearly parallel $\mathbf{G_2}$ structures.}  In
particular, it is also clear that when $\tau_1=0$ the associated
Riemannian metric is Einstein with positive scalar curvature.  This is
a well known fact of nearly parallel (type $1$) or weak holonomy $G_2$
manifolds, see~\cite{Gr2}.  The form of the curvature for a nearly
parallel $G_2$ structure was given as follows in~\cite{carrion,Cl}
\begin{lemma}
  The Riemannian curvature tensor of the metric associated to a nearly
  parallel $G_2$ structure has the form $\RC =
  W+\tfrac1{32}\tau_0^2r_g(g)$, where $W\in\mW_{77}$.
\end{lemma}

\noindent
{\bf Locally conformally parallel $\mathbf{G_2}$ structures.}  When on
the other hand $\tau_0=0$ one has a type $4$ structure.  Ricci and
scalar curvature in this case are non-trivial, in general.  Compact
manifolds with $G_2$ structure locally conformal to a parallel $G_2$
structure were described in \cite{math.DG/0509038,math.DG/0507179}. 

\noindent
{\bf Curvature of parallel $\mathbf{G_2}$ structures.}  When
$d\phi=0=d\Hodge\phi$, the only non-trivial component of the
Riemannian curvature is precisely $W_{77}$, see~\cite{MR0231313}.

\begin{remark}  
  An immediate consequence of the discussion above is that all $G_2$
  structures in the torsion class $1+4$ are generalized Einstein for
  $k=(4,-5)$.  Nearly parallel and parallel $G_2$ structures are
  generalized Einstein for all values of $k$.  These examples show how
  certain restrictions on the torsion gives restrictions in the form
  of the curvature.  As we shall see in section~\ref{sec:examples}
  drawing conclusions the other way around is not feasible.  Even when
  the form of the curvature is that of a parallel $G_2$ manifold,
  i.e., $\RC=W_{77}$, the only constraint on the torsion type is that
  of scalar-flatness.
\end{remark}

\section{Curvature of closed $G_2$ structures}
\label{sec:2}
A closed $G_2$ structure is by definition given by a closed $G_2$
three-form $\phi$.  Let us first examine the consequences of the
structure equations~\eqref{eq:16} in this case.  When $d\phi=0$ the
torsion $\tau$ has only one component $\tau=\tau_2$:
\begin{equation}
  \label{eq:32}
  d\Hodge\phi=\tau\bw\phi,
\end{equation}
whence $\delta^g\phi=\tau=\bar\xi$.  This is equivalent to the
equation $\bw_3(\xi)=0$ for the intrinsic torsion $\xi$ viewed as a
tensor in $\Lambda^2V^*\otimes V^*$.  This is just
the identity
\begin{equation}
  \label{eq:33}
  \xi_{ijk} + \xi_{jki} + \xi_{kij} = 0,
\end{equation}
for components of $\xi$ in an orthonormal frame.  Differentiating and
applying the Hodge star operator to the structure
equation~\eqref{eq:32} yields
\begin{equation}
  \label{eq:34}
  d\tau\bw\phi=0,\qquad\text{and} \qquad\delta^g\tau=0.
\end{equation}
The first equation is equivalent to $(d\tau)_7=0$ while the second may
be interpreted as $\LC\tau$ having no component in the $7$ dimensional
irreducible $\SO(7)$-submodule of $V^*\otimes\Lambda^2V^*$.  Using
Table~\ref{tab:g2-1} and Lemma~\ref{lem:alg} either of the
equations~\eqref{eq:34} is equivalent to the following Lemma.
\begin{lemma}\label{lem:dbar}
  Suppose $d\phi=0$.  Then $\Nt\tau\in V_{64}+V_{27}\subset
  V^*\otimes\Lambda^2_{14}$.\qed 
\end{lemma}
Equation~\eqref{eq:30} becomes
\begin{equation}
  \label{eq:35}
  d^{\Nt}\!\tau =  d\tau - \tfrac16\Hodge(\tau\bw\tau) -
  \tfrac16\abs{\tau}^2\phi,
\end{equation}
after rearranging.  Lemma~\ref{lem:dbar} shows that
$d^{\Nt}\!\tau\in\Lambda^3_{27}$.  This furnishes the observation:
\begin{equation}\label{eq:36}
  \tfrac13\Hodge d(\tau^3) = \inp{d\tau}{\Hodge(\tau\bw\tau)} = \langle {d^{\Nt}\!\tau},{\Hodge(\tau\bw\tau)_{27}}\rangle.
\end{equation}
The importance of this innocent looking set of relations is that the
left-hand side clearly integrates to zero on a compact manifold $M$.
We will use this in several places below.

\subsection{The Ricci curvature of a closed $G_2$ structure.}

Some of the next results have appeared in slightly different form and
for the special case of $k=(1,0)$ in~\cite{Br1} and~\cite{CI}.  The
main difference here concerns the interpretation of the results in
terms of the components of the covariant derivative $\Nt\tau$.
  
When the $G_2$ three-form is closed the formula for the generalized
Ricci curvature \eqref{eq:new1} may be written
\begin{align}\label{eq:37}
  \lambda_3\bigl(\Ric^k_0\bigr)
  &= -(k_1-4k_2)d^{\Nt}\!\tau +
  \tfrac13(k_1+5k_2)\Hodge(\tau\bw\tau)_{27}. 
\end{align}
This leads to
\begin{multline}
  \label{eq:38}
  \lVert{\Ric^k_0}\rVert^2 =
  \tfrac12(k_1-4k_2)^2\lvert{d^{\Nt}\!\tau}\rvert^2 +
  \tfrac1{21}(k_1+5k_2)^2\abs{\tau}^4\\ -
  \tfrac13(k_1+5k_2)(k_1-4k_2)\langle
  {d^{\Nt}\!\tau},{\Hodge(\tau\bw\tau)_{27}}\rangle.
\end{multline}
Here we have used equation~\eqref{eq:7} and the following relations,
valid for any two-form $\tau\in\Lambda^2_{14}$, see~\cite{Br1}
\begin{equation*}
  \Hodge(\tau\bw\tau\bw\phi) = -\abs{\tau}^2,\qquad \abs{\tau\bw\tau}^2 =
  \abs{\tau}^4 \qquad \text{and} \qquad \abs{(\tau\bw\tau)_{27}}^2 =
  \tfrac67\abs{\tau}^4.
\end{equation*}
All these may be verified by observing that, since the trivial
representation occurs in the symmetric powers $S^2\Lambda^2_{14}$ and
$S^4\Lambda^2_{14}$ with multiplicity one, the left- and right-hand
side of the equations must be proportional.  The constant of
proportionality is then found by evaluating on a sample element.
Alternatively, the last equation follows from the first two by
projecting $(\tau\bw\tau)_{27} = \tau\bw\tau +
\tfrac17\abs{\tau}^2\Hodge\phi$ and taking norms.

Recall the scalar curvature of the metric associated to a closed $G_2$
structure is $s_g=-\tfrac12\abs{\tau}^2$.  Integrate \eqref{eq:38}
over the compact manifold to get with the help of \eqref{eq:36}
\begin{proposition}\label{prop:2}
  Suppose $(M,\phi)$ is a compact $G_2$ manifold with $d\phi=0$.  Then
  \begin{equation}\label{bryp}
    \int_M \lVert{\Ric^k_0}\rVert^2 dV_g \geqslant \tfrac4{21} (k_1 +
    5k_2)^2 \int_M s_g^2 dV_g,  
  \end{equation}
  where equality holds if and only if $k_1=4k_2$ or $d^{\Nt}\!\tau=0$.
  In particular, if any one of $\Ric_0^g,~\Ric_0^\phi$ or $\Ric^\mW$
  is zero then $\tau=0$.
\end{proposition}
This includes the inequality established in ~\cite{Br1} for $k=(1,0)$.
Bryant studied equality in \eqref{bryp} for $k=(1,0)$ and referred to
the case, as that of an \emph{extremally pinched Ricci curvature}.
Such structures are characterized by
$$\lVert{\Ric^k_0}\rVert^2  =\tfrac4{21} s_g^2,
$$
which, in view of Proposition~\ref{prop:2} is equivalent to the
condition $d^{\Nt}\!\tau=0$. Bryant showed that the scalar curvature
is constant if the manifold is compact and gave a homogeneous compact
example.

Closely related statements are obtained by applying the
Cauchy-Schwartz inequality to the last summand of
equation~\eqref{eq:38}.  Here $M$ is not required to be compact.
\begin{corollary}\label{cor:3}
  Suppose $M$ is a manifold equipped with a closed $G_2$ structure
  $\phi$ with torsion $\tau$.  Let $k_1,k_2$ be real numbers and set
  $K_1:=\abs{k_1 - 4k_2}, ~K_2:=\abs{k_1 + 5k_2}$.  Then the
  inequality
  \begin{equation}\label{eq:40}
    \tfrac12\bigl(K_1\lvert{d^{\Nt}\!\tau}\rvert -
    \sqrt{\tfrac2{21}}K_2\abs{\tau}^2\bigr)^2 \leqslant
    \lVert{\Ric^k_0}\rVert^2 \leqslant
    \tfrac12\bigl(K_1\lvert{d^{\Nt}\!\tau}\rvert +  
    \sqrt{\tfrac2{21}}K_2\abs{\tau}^2\bigr)^2 
  \end{equation}
  holds everywhere in $M$.  If, in particular, $d^{\Nt}\!\tau=0$ then
  $\lVert{\Ric^k_0}\rVert^2 = \tfrac4{21}(k_1+5k_2)^2s^2_g$\qed
\end{corollary}

\begin{remark}
  Note that equation~\eqref{eq:37} and integration of \eqref{eq:38}
  over a compact $G_2$ manifold has an additional consequence: a
  compact $M$ with closed fundamental form is generalized Einstein for
  $k=(k_1,k_2)$ if and only if it is either parallel (in which case it
  is generalized Einstein for all $k$) or $k_1+5k_2=0$ and $\phi$ has
  extremally pinched Ricci curvature.  This result is implicitly
  contained in the discussion of Remark~14 in~\cite{Br1}.
\end{remark}
\subsection{The curvature component $\mathbf{W}_{64}$.}
\label{sec:pontry-class-clos}
Let $S\colon\so(n)\otimes\so(n)\to\so(n)\otimes\so(n)$ be defined by
$S(\alpha\otimes\beta)=\beta\otimes\alpha$. 
\begin{lemma}\label{prop:4}
  Suppose $\phi$ is a closed $G_2$ three-form with associated metric
  $g$.  Let $\tau$ be the torsion of $\phi$ and write $\Nt\tau_{64}$
  and $\Nt\tau_{27}$ for the respective components of $\Nt\tau$ in the
  ${64}$ and $7$ dimensional subspaces of $V^*\otimes\Lambda^2_{14}$.
  Then the component $W_{64}$ of the Riemannian curvature $\RC$
  satisfies
  \begin{equation*}
    W_{64}=S((\Nt\xi)_{64})+(\Nt\xi)_{64}=S(\Rt_{64})+\Rt_{64}.
  \end{equation*}
  Furthermore, $\norm{W_{64}}^2 = \tfrac13\lVert\Nt\tau_{64}\rVert^2$.
\end{lemma}
\begin{proof}
  Note that $W_{64}\in S^2(\so(7))$ so $S(W_{64})=W_{64}$.  Since
  $\xi\otimes\xi$ lies in a submodule of $S^2(V^*\otimes\lie
  g_2^{\perp})$ isomorphic to $S^2(\lie g_2)$, the tensor $(\xi^2)$
  does not contribute to $W_{64}$ and so
  $W_{64}=\Rt_{64}+(\Nt\xi)_{64}=S(\Rt_{64}+(\Nt\xi)_{64})$.  However,
  $\Rt_{64}\in\lie g_2^{\perp}\otimes\lie g_2$ while
  $(\Nt\xi)_{64}\in\lie g_2\otimes\lie g_2^{\perp}$, so
  \begin{equation*}
    \lie g_2^{\perp}\otimes\lie g_2\ni\Rt_{64}-S((\Nt\xi)_{64}) =
    (\Nt\xi)_{64}-S(\Rt_{64})\in\lie g_2\otimes\lie g_2^{\perp}. 
  \end{equation*}
  Therefore $\Rt_{64}=S((\Nt\xi)_{64})$ and
  $W_{64}=S((\Nt\xi)_{64})+(\Nt\xi)_{64}=S(\Rt_{64})+\Rt_{64}$.
  
  Let $e_i$ be a local orthonormal frame.  Using
  equation~\eqref{eq:17} and~\eqref{eq:19} we get
  \begin{align*}
    (\Nt\xi)_{ijkl}&=\Nt_i\xi_{jkl}-\Nt_j\xi_{ikl}
    =\tfrac16\bigl(\Nt_i\tau_{jp}-\Nt_j\tau_{ip})\phi_{pkl},
    =\tfrac16\bigl((d^{\Nt}\tau)_{ijp} - \Nt_p\tau_{ij})\phi_{pkl}\\
    &=\tfrac16\bigl((d^{\Nt}\tau)_{ijp} - (\Nt\tau_{27})_{pij} -
    (\Nt\tau_{64})_{pij})\phi_{pkl}.
  \end{align*}
  Here we make explicit use of the principle given by
  equation~\eqref{eq:13}.  Since $d^{\Nt}\tau\in\Lambda^3_{27}$ by
  Lemma~\ref{lem:dbar}, projection gives $ (\Nt\xi_{64})_{ijkl}= -
  \tfrac16(\Nt\tau_{64})_{pij}\phi_{pkl}$.  Take the tensor norm and
  use equation~\eqref{eq:8} to get
  \begin{align*}
    \lVert(\Nt\xi_{64})\rVert^2 &=
    \tfrac1{36}(\Nt\tau_{64})_{pij}\phi_{pkl}(\Nt\tau_{64})_{qij}\phi_{qkl},
    =\tfrac16(\Nt\tau_{64})_{pij}(\Nt\tau_{64})_{pij},
    =\tfrac16\lVert\Nt\tau_{64}\rVert^2.
  \end{align*}
  The final equation of the Lemma follows from this.
\end{proof}

\subsection{Closed fundamental three-forms with parallel torsion}

We briefly describe the only known example of a $G_2$
structure with closed fundamental form and extremally pinched Ricci
curvature on a compact manifold, due to Robert Bryant~\cite{Br1}.  Let
$G$ be the space of affine transformations of $\mathbb C^2$ preserving
the canonical complex volume form.  Then $\SU(2)$ is a subgroup in $G$
and $M = G/\SU(2)$ is a $7$-dimensional homogeneous space,
diffeomorphic to $\mathbb R^7$. It admits an invariant closed $G_2$
structure $\phi$ with extremally pinched Ricci curvature, as well as a
free and properly discontinuous action of a discrete subgroup
$\Gamma\subset G$ for which $\phi$ is invariant. So $\phi$ descends to
a closed $G_2$ structure $\tilde\phi$ with extremally pinched Ricci
curvature on the compact quotient $\tilde M:=\Gamma\setminus M$.

Andrew Swann made us aware of the following alternative description of
Bryant's example.  Note that $\SL(2,\mathbb C)=\SU(2).\Lie{Sol}_3$
where $\Lie{Sol}_3$ is the space of complex upper triangular $2\times
2$ matrices $\bigl(
\begin{smallmatrix}
  e^t & z\\ 0 & e^{-t}
\end{smallmatrix}\bigr)$, with $t$ real and $z$ a complex number
(by Iwasawa decomposition, or simply applying the Gram-Schmidt process
to the column vectors of $SL(2,\mathbb C)$).  This gives the
alternative description of $M$ as the Lie group $\Lie{Sol}_3 \rtimes
\mathbb C^2$.  Taking any basis of left-invariant one-forms $e=(e^i)$
on $M$ thought of as a Lie group gives a $G_2$ three-form $\phi$ by
requiring that $e$ is a $G_2$ adapted frame field and hence
multitude of left-invariant $G_2$ structures on $M$.  Below we shall
show that up to isometries there is only one closed three-form with
extremally pinched Ricci curvature among the many invariant
three-forms $M$ bears.

Even better, we can show that any manifold $M'$ with closed $G_2$
structure $\phi'$ such that the torsion $\tau$ is parallel with
respect to $\Nt$, is locally isometric to $(M,\phi)$.
Lemma~\ref{lem:dbar} and Proposition~\ref{prop:4} give
\begin{lemma}
  A closed $G_2$ form $\phi$ has extremally pinched Ricci curvature if
  and only if its torsion $\tau$ satisfies $(\Nt\tau)_{27}=0$.  The
  torsion of a closed $G_2$ form $\phi$ is $\Nt$-parallel if and only
  if $\phi$ has extremally pinched Ricci curvature and the component
  $W_{64}$ of the Riemannian curvature is zero.
\end{lemma}

\subsection*{Proof of Theorem~\ref{thm:par-tors}a.}
If $\Nt\tau=0$ then $(M,\phi)$ has extremally pinched Ricci curvature
and the norm of the torsion is constant.  From this it follows that
$\tau$ has constant rank $4$, see Remark~13 of~\cite{Br1}.  The
stabilizer algebra of $\tau$ in $\lie g_2$ is then isomorphic to
$\un(2)$ and the holonomy algebra of $\Nt$ accordingly reduces to a
$\un(2)$ in $\lie g_2$.  There are two such, distinguished as follows.
The subalgebras of $\lie g_2$ of dimension $4$ are both contained in
the maximal $\so(4)$ and both isomorphic to $\un(2)$.  The action of
$\so(4)\subset\lie g_2$ on $\mathbb R^7$ may be written in terms of
the standard representations $V_+,V_-$ of $\so(4)=\su_+(2)+\su_-(2)$
as $\mathbb R^7=S^2V_++ V_+\otimes V_-$ (modulo complexifications).
Take a basis $e_1,e_2,e_3,e_4$ of $V_+\otimes V_-\cong\mathbb R^4$ and
$e_5,e_6,e_7$ of $S^2V_+$.  Note that under the action of $\so(4)$ the
space of two-forms decomposes as
  \begin{gather*}
    \Lambda^2=2S^2V_++(S^3V_++V_+)\otimes V_-+S^2V_-,\quad\text{whence}\quad
    \Lambda^2_{14} = S^2V_+ + S^2V_- + S^3V_+\otimes V_-.
  \end{gather*}
  Explicit bases of the subspaces $S^2V_\pm$ are $h :=
  e^{12}+e^{34}-2e^{56}, x:= e^{13}+e^{42}-2e^{67}, y:=
  -(e^{14}+e^{23}-2e^{57})$ for $S^2V_+$ and $t:= e^{12}-e^{34}$ and
  so on for $S^2V_-$. Note that elements of $S^2V_+$ always have rank
  $6$.  Therefore the stabilizer of a rank $4$ element of $\lie g_2$
  is the algebra $\un_+(2)=\mathbb
  R_-+\su_+(2)\subset\so(4)\subset\lie g_2$.  This acts on $\mathbb
  R^7$ by $S^2V_++(L+\bar L)V_+$ and on $\Lambda^2$ by
  \begin{gather*}
    2S^2V_++(L^2+\mathbb R+\bar L^2)+(L+\bar L)(S^3V_++V_+).
  \end{gather*}
  Here $L$ is the standard complex representation of $\mathbb
  R=\un(1)$ on $\mathbb C$.  The Lie algebra $\un_+(2)$ has two
  non-trivial ideals isomorphic to $\un(1)$ and $\su(2)$.  The induced
  representations on $\mathbb R^7$ are $\mathbb R^3\oplus 2(L+\bar L)$
  and $S^2V_+ +2V_+$, respectively.  Neither is an irreducible
  holonomy representation and so $\mK(\un_+(2))$ is trivial,
  see~\cite{MR2114426} for further details.  By assumption $\Nt\xi=0$
  and so~\ref{thm:curv} may be applied to conclude that $M$ with its
  induced metric $g$ is locally isometric to some homogenous space.
  
  Which homogeneous space may be determined by the Nomizu
  construction, see~\cite{MR0059050}.  Given $\xi$ and $\Rt$ a Lie
  bracket is constructed on $\lie g:=\un_+(2)\oplus\mathbb R^7$ by
  setting $[A+x,B+y] : = ([A,B] -\Rt_{x,y}) + (Ay-Bx -\xi_xy +
  \xi_yx)$.  Now the universal covering space of $M$ is the
  homogeneous space $G/H$ where $G$ is the connected and simply
  connected Lie group with Lie algebra $\lie g$ and $H$ is the
  connected subgroup of $G$ with algebra $\un_+(2)$.  Note that $H$ is
  closed in $G$ by virtue of the completeness of $g$.  The induced
  action of $H$ on $\mathbb R^7$ is as the stabilizer group of the
  pair $(\phi,\tau)$, whence $H\cong\Un(2)\subset G_2$.  To calculate $G$
  we fix $\tau:= 6t$.  The elements $t,h,x,y\in\lie
  g_2\cong\Lambda^2_{14}$ defined above can then be taken as
  generators for $\un_+(2)$.  The intrinsic torsion $\xi$ is obtained
  through~\eqref{eq:3}.  The curvature operator $\Rt$ is calculated by
  noting that $\RC=\Rt+(\xi^2)$ (c.f.~\eqref{eq:20}) must satisfy the
  Bianchi identity whence $\bian\Rt=-\bian((\xi^2))$ .  Explicitly,
  the two subspaces of $\Lambda^2$ isomorphic to $S^2V_+$ are
  generated by
  \begin{gather*}
    a_1=e^{56},\quad b_1 = e^{67},\quad c_1 = e^{75},\quad
    \text{and}\quad a_2=e^{12}+e^{34},\quad b_2 = e^{13}+e^{42},\quad
    c_2 = -(e^{14}+e^{23}).
  \end{gather*}
  By $\un_+(2)$ invariance of $\Rt$ it then follows that constants $A,B$
  and $C$ exist such that
  \begin{multline}\label{rtcu}
    \Rt=At\otimes t + B(a_1\otimes h + b_1\otimes x +
    c_1\otimes y) + C(a_2\otimes h +  b_2 \otimes x + c_2\otimes y).
  \end{multline}
  On the other hand, using \eqref{eq:33} several times gives
  $\bian((\xi^2))_{ijkl} = \xi_{pij}\xi_{lkp} + \xi_{pjk}\xi_{lip} +
  \xi_{pki}\xi_{ljp}$.  Therefore $\bian(\Rt)_{ijkl} = -
  (\xi_{pij}\xi_{lkp}+\xi_{pjk}\xi_{lip}+\xi_{pki}\xi_{ljp})$.
  Evaluating left- and right- hand side now gives $b(\Rt)_{1234} = -A
  + 3C = 0,\quad b(\Rt)_{1256} = -2C = 0,\quad b(\Rt)_{5612} = B =
  -2$, whence $\Rt= -2(a_2\otimes h + b_2\otimes x + c_2\otimes y)$.
  As a map $\so(7)\to\un_+(2)$ the range of $\Rt$ is then $\su_+(2)$.
  Note that this allows to replace $\un_+(2)$ by $\su_+(2)$ in the
  construction of $\lie g$.  It is now a trivial exercise to write
  down the Lie brackets on $\su_+(2)+\mathbb R^7$ and verify that,
  indeed, $\lie g\cong\lie{sl}(2,\mathbb C)\rtimes \mathbb C^2$.
  \qed

\begin{remark}
  Many examples of homogeneous spaces (with compact quotients)
  admitting invariant closed $G_2$-structures are known,
  see~\cite{MR906398, MR887944}.  Generically these do not admit a
  closed three-form with extremally pinched Ricci curvature, by
  Theorem~\ref{thm:par-tors}.  Due to~\eqref{eq:35} in fact, the
  condition of having extremally pinched Ricci curvature may be seen
  as a non-degeneracy requirement on the derivative $d\tau$, which is,
  again generally speaking, not satisfied on examples built from
  nilpotent and solvable Lie algebras.
  
  The space $G/\SU(2)$ with its invariant closed three-form with
  extremally pinched Ricci curvature may be described in several
  different ways as we have already seen.  A very concrete way is
  provided by the rank one solvable extension of the complex
  Heisenberg group given by the structure equations
  \begin{equation*}
    \begin{array}{l}
      de^1=-e^{17}-2(e^{36}+e^{45}),\qquad
      de^2=-e^{27}-2(e^{35}+e^{64}),\\
      de^3=e^{37},\quad de^4=e^{47},\quad de^5=-2e^{57},\quad
      de^6=-2e^{67},\quad de^7=0.
    \end{array}
  \end{equation*}
  Here the closed $\phi$ with extremally pinched Ricci curvature may
  taken to be the left-invariant three-form given by
  equation~\eqref{eq:2}. 
\end{remark}

\section{Topology of  closed $G_2$ structures. Proof of  Theorem~\ref{thm:par-tors}b and Theorem~\ref{main1}b}
Closed $G_2$ structures are distinguished also by certain naturally
associated topological data.
\begin{theorem}\label{thm:pontrj}
  Suppose $M$ is a compact $7$ dimensional manifold with a closed
  fundamental three-form $\phi$.  Then
  \begin{equation}
    \label{eq:43}
    \inp{p_1(M)\cup[\phi]}{[M]}=-\frac1{8\pi^2}\int_M
    \left\{\norm{W_{77}}^2-\tfrac12\norm{W_{64}}^2 -
      \tfrac97\norm{\Ric^g_0}^2 + \tfrac{45}{28^2}
      s_g^2\right\} dV_g,
  \end{equation}
  where $p_1(M)$ is the first Pontrjagin class of $M$.
\end{theorem}
This generalizes Proposition~10.2.7. of~\cite{Joyce:book}.
\begin{proof}
  We shall be working in a local orthonormal frame $e_i$.  First,
  equations~\eqref{eq:8} and~\eqref{eq:33} yield
  \begin{multline*}
    \RC_{ijab} = \Rt_{ijab} + \tfrac16(\Nt_i\tau_{jp} -
    \Nt_j\tau_{ip})\phi_{pab} +
    \tfrac1{36}\tau_{pq}\tau_{pr}\phi_{qij}\phi_{rab}
    +\tfrac1{36}(\tau_{ia}\tau_{jb} - \tau_{ib}\tau_{ja})
    -\tfrac1{18}\tau_{ip}\tau_{jq}\phi_{pqab}.  
  \end{multline*}
  This relates the Riemannian curvature $\RC$ of the metric $g$
  associated to $\phi$ with the curvature $\Rt$ of the canonical
  connection $\Nt$ of $\phi$ with torsion expressed in terms of $\tau$
  rather than the intrinsic torsion $\xi$.
  Note that from Chern-Weil theory $p_1(M)$ may be represented by the
  $4$-form
  \begin{align*}
    \tfrac1{8\pi^2}\tr(\RC\bw\RC) =
    \tfrac1{16\pi^2}\RC_{ijab}\RC_{klab}e^{ijkl}.
  \end{align*}
  Now
  \begin{multline*}
    8\pi^2\inp{p_1(M)\cup[\phi]}{[M]} = \int_M
    \tr(\RC\bw\RC)\bw\phi
    = \int_M\inp{\tr(\RC\bw\RC)}{\Hodge\phi}dV_g,\\
    = \tfrac1{2}\int_M\RC_{abij}\RC_{cdij}\phi_{abcd}dV_g
    = -\int_M\norm{\RC}^2dV_g + \tfrac12\int_M
    (\RC_{ijab}\phi_{abt})(\RC_{ijcd}\phi_{cdt})dV_g. 
  \end{multline*}
  The contraction $\RC_{ijab}\phi_{abt}$ gives
  \begin{align*}
    \RC_{ijab}\phi_{abt} &= \tfrac16(\Nt_i\tau_{jp} -
    \Nt_j\tau_{ip})\phi_{pab}\phi_{abt} +
    \tfrac1{36}\tau_{pq}\tau_{pr}\phi_{qij}\phi_{rab}\phi_{abt} \\ &
    \quad + \tfrac1{36}(\tau_{ia}\tau_{jb} -
    \tau_{ib}\tau_{ja})\phi_{abt}
    -\tfrac1{18}\tau_{ip}\tau_{jq}\phi_{pqab}\phi_{abt}\\
    &= \Nt_i\tau_{jt} - \Nt_j\tau_{it} +
    \tfrac1{6}\tau_{pq}\tau_{pt}\phi_{qij}
    -\tfrac16\tau_{ip}\tau_{jq}\phi_{pqt}.
  \end{align*}
  where the identities~\eqref{eq:8}-~\eqref{eq:12} are applied.  By
  definition $\bw_3(\Nt\tau) =: d^{\Nt}\!\tau$ so
  \begin{equation}\label{eq:44}
    \RC_{ijab}\phi_{abt} = \bigl((d^{\Nt}\!\tau)_{ijt} - \Nt_t\tau_{ij}\bigr) +
    \tfrac1{6}(\tau_{pq}\tau_{pt}\phi_{qij} -
    \tau_{ip}\tau_{jq}\phi_{pqt})
  \end{equation}
  The evaluation of the integrand $(\RC_{ijab}\phi_{abt})
  (\RC_{ijcd}\phi_{cdt})$ is now reduced to evaluating $9$ different
  contractions.  Some of these are easy, for instance
  \begin{gather*}
    [(d^{\Nt}\!\tau)_{ijt} - \Nt_t\tau_{ij}][(d^{\Nt}\!\tau)_{ijt} -
    \Nt_t\tau_{ij}] = 2\lvert d^{\Nt}\!\tau\rvert^2 +
    \lVert\Nt\tau\rVert^2,\\
    \tau_{pq}\tau_{pt}\phi_{qij}\tau_{rs}\tau_{rt}\phi_{sij} =
    6\tau_{pq}\tau_{pt}\tau_{rq}\tau_{rt} =
    6\abs{\tau}^4.
  \end{gather*}
  The last equality is obtained by the standard method: observe that
  the right hand side is a fourth order homogeneous $G_2$ invariant
  polynomial in $\tau\in\Lambda^2_{14}$.  So up to scale it must equal
  $\abs{\tau}^4$.  The constant of proportionality is found by
  evaluating on a test element.  In the same way one obtains \(
  \tau_{ip}\tau_{jq}\phi_{pqt}\tau_{ir}\tau_{js}\phi_{rst} =
  3\abs{\tau}^4 \) and \( \tau_{pq}\tau_{pt}\phi_{qij}
  \tau_{ir}\tau_{js}\phi_{rst} = 0 \).  To evaluate the remaining
  terms first note that
  \begin{align*}
    (\Hodge(\tau\bw\tau)+\abs{\tau}^2\phi)_{ijt} &=
    \tau_{pq}(\tau_{pt}\phi_{qij} + \tau_{pi}\phi_{qjt} +
    \tau_{pj}\phi_{qti})\\
    &= 2(\tau_{ip}\tau_{jq}\phi_{pqt} + \tau_{jp}\tau_{tq}\phi_{pqi} +
    \tau_{tp}\tau_{iq}\phi_{pqj}).
  \end{align*}
  Lemma~\ref{lem:dbar} and relation~\eqref{eq:36} then imply that \(
  (d^{\Nt}\tau)_{ijt}\tau_{pq}\tau_{pt}\phi_{qij} = \tfrac23\Hodge
  d(\tau^3), \) and \( (d^{\Nt}\tau)_{ijt}\tau_{ip}\tau_{jq}\phi_{pqt}
  = \tfrac13\Hodge d(\tau^3).  \) The contraction
  $\Nt_k\tau_{ij}\phi_{ijl}$ is zero for all $k$ and $l$ as $\Nt\tau
  \in V^*\otimes\Lambda^2_{14}$.  So $\Nt_t\tau_{ij} \tau_{pq}
  \tau_{pt} \phi_{qij} = 0$.  Using this identity and $\bw_3(\xi)=0$
  repeatedly, gives
  \begin{align*}
    \Nt_t\tau_{ij}\tau_{ip}\tau_{jq}\phi_{pqt} &=
    -\Nt_t\tau_{ij}\tau_{ip}(\tau_{pq}\phi_{tqj} + \tau_{tq}\phi_{jqp})
    = \Nt_t\tau_{ij}\tau_{pq}\tau_{pi}\phi_{qjt} +
    \Nt_t\tau_{ij}\tau_{tq}(\tau_{jp}\phi_{qip} +
    \tau_{qp}\phi_{ijp})\\
    &= \Nt_t\tau_{ij}\tau_{pq}\tau_{pi}\phi_{qjt} -
    \Nt_t\tau_{ij}\tau_{jp}(\tau_{iq}\phi_{qpt} +
    \tau_{pq}\phi_{qti})\\
    &= \Nt_t\tau_{ij}\tau_{pq}(\tau_{pi}\phi_{qjt} +
    \tau_{pj}\phi_{qti} + \tau_{pt}\phi_{qij}) +
    \Nt_t\tau_{ij}\tau_{ip}\tau_{jq}\phi_{pqt}.
  \end{align*}
  Thus, $\Nt_t\tau_{ij}\tau_{ip}\tau_{jq}\phi_{pqt} = \tfrac13\Hodge
  d(\tau^3)$. The net result is
  \begin{equation*}
    (\RC_{ijab}\phi_{abt})(\RC_{ijcd}\phi_{cdt}) = 2\lvert
    d^{\Nt}\!\tau\rvert^2 + \lVert\Nt\tau\rVert^2 +
    \tfrac14\abs{\tau}^4 + \tfrac29\Hodge d(\tau^3). 
  \end{equation*}
  Using equation~\eqref{eq:38} with $k=(1,0)$ and Lemmas~\ref{lem:alg}
  and~\ref{prop:4}, this may be reformulated as
  \begin{equation}\label{nca}
    (\RC_{ijab}\phi_{abt})(\RC_{ijcd}\phi_{cdt}) = 3\norm{W_{64}}^2 +
    \tfrac{40}7\lVert\Ric^g_0\rVert^2 - \tfrac{13}{147}s_g^2
    + \tfrac{34}{63}\Hodge d(\tau^3).
  \end{equation}
  Take $k=(1,0)$ and $k=(4,-5)$ in \eqref{eq:38} and make a simple
  substitution using \eqref{eq:36}.  Now apply Stokes' Theorem to the
  obtained equalities.  This gives
  \begin{gather}
    \label{eq:39}
    \int_M\left(36\norm{\Ric_0^g}^2 - 25\lVert\Ric^{\mW}\rVert^2 -
      \tfrac{45}{28}s_g^2\right)dV_g = 0.  
  \end{gather}
  Integrating \eqref{nca}, applying Theorem~\ref{thm:curv-decomp} and
  \eqref{eq:39} we finally obtain \eqref{eq:43}.
\end{proof}

\begin{proposition}\label{prop:3}%
  Suppose $M$ is a compact $7$-dimensional manifold equipped with a
  closed $G_2$ structure $\phi$.  Then
  \begin{enumerate}[$(a)$]
  \item \begin{equation*}
    \inp{p_1(M)\cup[\phi]}{[M]} \geqslant -\frac1{8\pi^2}\int_M
    \left\{\norm{W_{77}}^2 - \tfrac12\norm{W_{64}}^2 - \tfrac{3}{16}
    s_g^2\right\} dV_g.
  \end{equation*}
  Equality holds if and only if $\phi$ has extremally pinched Ricci
  curvature.\label{item:2}
\item \begin{equation*} \inp{p_1(M)\cup[\phi]}{[M]} \geqslant
    -\frac1{8\pi^2}\int_M \left\{\norm{W_{77}}^2 - \tfrac{3}{16}
      s_g^2\right\} dV_g.
  \end{equation*}
  Equality holds if and only if $\Nt\tau=0$.
\item \begin{equation*} \inp{p_1(M)\cup[\phi]}{[M]} \geqslant
    -\frac1{8\pi^2}\int_M \norm{W_{77}}^2 dV_g,
  \end{equation*}
  where equality holds if and only if $\tau=0$.\label{item:3}
  \end{enumerate}
\end{proposition}
\begin{proof}
  Using \eqref{eq:38} in  formula~\eqref{eq:43} gives
    \begin{equation*}
      \langle p_1(M)\cup[\phi],[M]\rangle =
      -\frac1{8\pi^2}\int_M\left\{\norm{W_{77}}^2 -
        \tfrac12\lVert{W_{64}}\rVert^2 - \tfrac9{14}\lvert
        d^{\Nt}\tau \rvert^2 -\tfrac3{16}s_g^2\right\} dV_g.
    \end{equation*}
    This gives the first statement.  The second follows by
    Lemma~\ref{prop:4}.  The third is now a consequence of the scalar
    curvature formula~\eqref{eq:29} for a closed $G_2$ structure.
\end{proof}

\begin{proposition}
  Suppose $M$ is compact, equipped with a closed $G_2$ three-form
  $\phi$.  Suppose furthermore that the component $W_{77}$ of the
  Riemannian curvature is zero.  Then
  \begin{equation*}
    \inp{p_1(M)\cup[\phi]}{[M]} \geqslant 0 
  \end{equation*}
  and equality holds if and only if $\phi$ is parallel and the
  associated metric is flat.  
  
  If, furthermore, the cohomology class $[\phi]$ contains a parallel
  $G_2$ three-form $\phi'$ then equality is attained and $\phi=\phi'$.
\end{proposition}
\begin{proof}
  When $W_{77}=0$, Proposition~\ref{prop:3}~\ref{item:3} gives
  $\inp{p_1(M)\cup[\phi]}{[M]} \geqslant 0$
  with equality if and only if $\tau=0$.  If $\tau=0$ then
  $W_{64}=0=W_{27}$ and $\Ric^g=0$ (see section~\ref{sec:1+4}).  This
  gives the first statement of the Corollary.
  
  Suppose now that $\phi$ is closed with $W_{77}=0$ and that a
  parallel $G_2$ three-form $\phi'$ exists in $[\phi]$.  Then
  $0\geqslant
  \inp{p_1(M)\cup[\phi']}{[M]}=\inp{p_1(M)\cup[\phi]}{[M]}\geqslant
  0$.  Therefore the torsion of $\phi$ is zero and $\phi=\phi'$.
\end{proof}
In particular, suppose $\phi$ is a closed $G_2$ structure on a compact
manifold $M$.  If $\beta$ is a two-form on $M$ that is not closed such
that $\tilde\phi = \phi+d\beta$ is a $G_2$ three-form then the metric
$\tilde g$ associated to $\tilde\phi$ has $W_{77}^{\tilde g}\not=0$.
This is independent of the value of $W_{77}$ for $\phi$.

\subsection{Proofs of main theorems}
\begin{proof}[Proof of Theorem~\ref{thm:par-tors}b]
  The stated result is now an easy consequence of
  Theorem~\ref{thm:par-tors}a and Corollary
\end{proof}
\begin{proof}[Proof of Theorem~\ref{main1}]
  Suppose $W_{27}=0$.  Then Theorem~\ref{thm:curv-decomp} and
  Proposition~\ref{prop:2} yield $\phi$ is parallel.  This gives
    the first statement of Theorem~\ref{main1}.  The second
    statement is a consequence of the following two Lemmas.
\begin{lemma}[Integral Weitzenb\"ock Formula for $\tau$]\label{Weitz} 
  Let $\tau$ be the torsion of a closed $G_2$ structure $\phi$ on a
  compact manifold $M$.  Then
  \begin{equation}\label{eq:45}
     \int_M 4\RC(\tau,\tau)dV_g =   \int_M (2\lvert{d^{\Nt}\!\tau}\rvert^2 - 
    \lVert{\Nt\tau}\rVert^2 + \tfrac16\abs{\tau}^4)dV_g. 
  \end{equation}
\end{lemma}
\begin{proof}
  This may also be found in~\cite{CI}.  However, as conventions are
  different we give an outline.  Take the covariant derivative of the
  one-form $\Theta$ with components $\RC_{ijkl}\phi_{pij}\tau_{kl}$
  (expressed in terms of a local orthonormal frame) with respect to
  the Levi-Civita connection and contract with the metric.  The result
  is clearly a divergence and so, by Stokes' Theorem integrates to
  zero on the compact $M$.  But we also have $ \delta^g\Theta =
  -(\LC_p\RC_{ijkl})\phi_{pij}\tau_{kl} -
  \RC_{ijkl}\LC_p\phi_{pij}\tau_{kl} -
  \RC_{ijkl}\phi_{pij}\LC_p\tau_{kl}$.  Writing $\delta^g\phi=\tau$
  and applying the second Bianchi identity to the factor
  $\LC_p\RC_{ijkl}\phi_{pij}$ gives $\delta^g\Theta =
  \RC_{ijkl}\tau_{ij}\tau_{kl} - \RC_{ijkl}\phi_{pij}\LC_p\tau_{kl}$.
  The proof is finished by expressing $\LC\tau$ as a covariant
  derivative with respect to the canonical connection and then using
  the identities needed to prove equation~\eqref{eq:43}.
\end{proof}

The integrand $\RC(\tau,\tau)$ may be evaluated using
Theorem~\ref{thm:curv-decomp}, while the norms of the derivatives and
$\tau$ may be expressed in terms of curvature components using
formula~\eqref{eq:38} and Lemma~\ref{prop:4}.  This gives
\begin{lemma}\label{prop:5}
  Let $(M,\phi)$ be a compact $G_2$ manifold $M$ with closed $G_2$
  structure $\phi$ and torsion $\tau$.  Then the following integral
  identity holds
  \begin{align*}
  \int_M(4 W_{77}(\tau,\tau) +3\norm{W_{64}}^2)dV_g & =
  \int_M\left(\tfrac{16}7\lVert\Ric^g_0\rVert^2 +
    \tfrac{127}{196}s^2_g\right)dV_g.
\end{align*}
\end{lemma}
\begin{proof}
  A computation in an orthogonal frame shows that $ r_g(h)(\tau,\tau)=
  -h_{pq}\tau_{pr}\tau_{qr}$, $r_g(g)(\tau,\tau) = -2\abs{\tau}^2 =
  4s_g$ and $ r_\phi(h) = -\tfrac23h_{pq}\tau_{pr}\tau_{qr} +
  \tfrac13\tr_g(h)\abs{\tau}^2$ for any symmetric two-tensor $h$.
  This gives
  \begin{equation}\label{eq:46}
    4\RC(\tau,\tau) = 4W_{77}(\tau,\tau) - \tfrac14
    \bigl(\Ric_0^{(3,1/4)})_{pq}\tau_{pr}\tau_{qr} +
    \tfrac4{21}s^2_g. 
  \end{equation}
  The contraction $\tau_{pr}\tau_{qr}$ gives the components of the
  tensor $\tau\otimes_g\tau := (e_r\hook\tau)\otimes (e_r\hook\tau) =
  \tfrac12\tau_{pr}\tau_{qr}e^p\odot e^q$.  It is easy to check that $
  \lambda_3(\tau\otimes_g\tau) = \Hodge(\tau\bw\tau) +
  \abs{\tau}^2\phi$, and $ \lambda_3\bigl(\Ric_0^{(3,1/4)}\bigr) = -2
  d^{\Nt}\!\tau + \tfrac{17}{12}\Hodge(\tau\bw\tau)_{27}$ follows
  from~\eqref{eq:37}.  Relating inner products of two-tensors and
  three-forms via~\eqref{eq:7} and using equation~\eqref{eq:36}
  gives
  \begin{align}
    \label{eq:47}
    \bigl(\Ric_0^{(3,1/4)})_{pq}\tau_{pr}\tau_{qr} &= \tfrac12\inp{-2
      d^{\Nt}\!\tau +
      \tfrac{17}{12}\Hodge(\tau\bw\tau)_{27}}{\Hodge(\tau\bw\tau)_{27}}
    = -\tfrac1{3} \Hodge d(\tau^3) - \tfrac{17}{7}s_g^2,
  \end{align}
  whence $ 4\RC(\tau,\tau) = 4W_{77}(\tau,\tau) - \tfrac{5}{12}s^2_g +
  \tfrac1{12} \Hodge d(\tau^3)$, and clearly
  \begin{equation}\label{eq:48}
    \int_M 4\RC(\tau,\tau)dV_g = \int_M (4W_{77}(\tau,\tau) -
    \tfrac{5}{12}s^2_g)dV_g. 
  \end{equation}
  On the other hand,
  \begin{align*}
    2\lvert{d^{\Nt}\!\tau}\rvert^2 - \lVert{\Nt\tau}\rVert^2 +
    \tfrac16\abs{\tau}^4 &= - \lVert{(\Nt\tau)_{64}}\rVert^2 +
    (2-\tfrac67)\bigr(2\lVert\Ric^g_0\rVert^2 -\tfrac8{21}s^2_g +
    \tfrac29\Hodge d(\tau^3)\bigl) + \tfrac23
    s_g^2\\
    &= - 3 \norm{W_{64}}^2 + \tfrac{16}7\lVert\Ric^g_0\rVert^2
    +\tfrac{2\cdot17}{3\cdot 49}s_g^2 + \tfrac{16}{63}d(\tau^3).
  \end{align*}
  This gives an alternative expression for the Weitzenb\"ock formula
  \begin{align*}
    \int_M 4\RC(\tau,\tau) & = \int_M (- 3 \norm{W_{64}}^2 +
    \tfrac{16}7\lVert\Ric^g_0\rVert^2 +\tfrac{2\cdot17}{3\cdot
      49}s_g^2).
  \end{align*}
  We compare this to equation~\eqref{eq:48} and rearrange to get
  \begin{align*}
    \int_M(4 W_{77}(\tau,\tau) +3\norm{W_{64}}^2)dV_g & =
    \int_M\left(\tfrac{16}7\lVert\Ric^g_0\rVert^2 +
      \tfrac{127}{196}s^2_g\right)dV_g.
  \end{align*}
\end{proof}
Eventually, suppose $W_{77}=W_{64}=0$ and apply Lemma~\ref{prop:5} to complete the
proof of Theorem~\ref{main1}.
\end{proof}
\begin{corollary}
  Let $M$ be a compact $G_2$ manifold with closed $G_2$ structure
  $\phi$ and torsion $\tau$.  Then
  \begin{equation*}
    \int_M r_g(\Ric^g)(\tau,\tau)dV_g = 0.
  \end{equation*}
\end{corollary}
\begin{proof}
  A computation similar to~\eqref{eq:47} of the proof of
  Lemma~\ref{prop:5} shows that $\tfrac14 r_g(\Ric^g)(\tau,\tau) =
  \Ric^g_{pq}\tau_{pr}\tau_{qr} = -\tfrac16\Hodge d(\tau^3)$.
  Integrate the latter over the compact $M$ to get the assertion.
\end{proof}
\begin{remark}
  The name we have given Lemma~\ref{Weitz} stems from the equivalence
  of equation~\eqref{eq:45} with the usual integral Weitzenb\"ock
  formula $\int_M \abs{d\tau}^2 dV_g = \int_M
  \tfrac12\left(\norm{\LC\tau}^2 + 2\RC(\tau,\tau)\right)dV_g$,
  see~\cite{CI}.  This is special instance of an example considered
  in~\cite{math.DG/0702031}, namely the Weitzenb\"ock formula for two
  forms in $\Lambda^2_{14}$.
\end{remark}

\section{Examples}
\label{sec:examples}

All instances considered are locally of the form $M=I\times M^*$ where
$M^*$ is a $6$-dimensional manifold carrying a one-parameter family of
$\SU(3)$ structures given by $(\omega_t,\psi^+_t)_{t\in I}$ and $I$ is
a real interval.  This gives $G_2$ structures
$\phi:=dt\bw\omega_t+\psi^+_t$ with associated metric $g=dt^2+g_t$,
where $g_t$ is the metric on $M^*$ determined by $\omega_t,~\psi^+_t$.
Strictly speaking, pull-backs such as $g=\pi_1^*(dt^2)+\pi_2^*(g_t)$
where $\pi_i$ is projection to the $i$'th factor ought to be included,
but to simplify notation we set $dt:=\pi^*_1(dt)$ et cetera.

The examples all admit parallel or nearly parallel $G_2$ structures.
Some are flat, others have constant sectional curvature.  In
particular, the curvature can be written $W_{77}+S$, and sometimes
$W_{77}=0$ (constant sectional curvature), whilst in other (parallel)
cases $S=0$.  The compatible three-forms, however, appear to range
over pretty much any torsion type not excluded by the value of the
scalar curvature.

\subsection{Warped products}

Let $M^*$ be an $n-1$ dimensional manifold with metric $g^*$.  Write
${\RC}^*$, $\Ric^*$ and $(n-1)(n-2)\rho^*=s^*$ for the Riemannian, Ricci
and scalar curvatures of $g^*$.  Set $M:=I\times M^*$ with warped
product metric $g:=dt^2+f^2 g^*$ and curvatures
$\RC,~\Ric,~n(n-1)\rho=s$.  An elementary calculation shows that the
curvatures of $g$ are related to (the pull-backs relative to $M\to M^*$
of) those of $g^*$ via
\begin{equation}
  \label{eq:49}
  \RC=f^2{\RC}^*-\tfrac12(ff')^2 g^*\ovee g^* - ff'' dt^2\ovee g^*,
\end{equation}
It follows directly from equation~\eqref{eq:49} that
\begin{lemma}\cite{Besse:Einstein}
  The warped product metric $g=dt^2+f^2g^*$ is Einstein if and only if
  $g^*$ is Einstein and $ (f')^2+\rho f^2=\rho^*$. If $g$ is Einstein,
  then \( f''+\rho f=0, \) and $ \RC = f^2W^*+\tfrac12\rho g\ovee g, $
  where $\RC$ is the Riemannian curvature of $g$ and $W^*$ is the Weyl
  curvature of $W^*$ considered as $(4,0)$ tensors.
\end{lemma}

\subsection{Curvature of nearly K\"ahler $3$-folds}

Let $M^*$ be a $6$-dimensional manifold equipped with an $\SU(3)$
structure, i.e., with data $g^*,J,\omega,\psi^+,\psi^-$ where $g^*$ is
a Riemannian metric, $J$ is an almost complex structure, $\omega$ is a
non-degenerate form, and $\psi^++i\psi^-$ is a complex volume form (of
type $(3,0)$).  The normalization conditions
\begin{gather}
  \label{eq:51}
  \omega= g\circ J,\qquad 2\omega^3=3\psi^+\bw\psi^-,
\end{gather}
may be imposed to give relations $\Hodge\psi^+=\psi^-,\quad
J\psi^+=-\psi^-$, and so on.  The structure equations of a nearly
K\"ahler manifold are
\begin{equation}\label{eq:52}
  d\omega=3\sigma\psi^+,\qquad d\psi^-=-2\sigma\omega^2,
\end{equation}
where $\sigma$ is a positive constant related to scalar curvature
$s^*$ and normalized scalar curvature $\rho^*$ through $\sigma^2 =
s^*/30 = \rho^*$.  The Riemannian curvature tensor takes the form
similar to the one of a $G_2$ structure of type $1$ 
\begin{gather*}
  {\RC}^* =  K^* + \tfrac12\rho^*(g^*\ovee g^*),
\end{gather*}
where $K^*(\omega)=0$, $K^*(\iota_X\psi^{\pm})=0$ for all $X\in
\Gamma(M^*)$, see~\cite{carrion}.  
 
\subsection{Warped Products over nearly K\"ahler $3$-folds}
\label{sec:warped-products-over}
Let $I\subset \mathbb R$ be an open interval and set $M=M^*\times I$
where $M^*$ has an $\SU(3)$ structure $(g^*,J,\omega,\psi^+,\psi^-)$.
Define
\begin{gather*}
  \omega_t=f^2\omega,\quad
  \psi_t^+=f^3(\cos\theta\psi^+-\sin\theta\psi^-),\quad
  \psi_t^-=f^3(\sin\theta\psi^++\cos\theta\psi^-),
\end{gather*}
for smooth functions $f,\theta\colon I\to\mathbb R$ with $f>0$.  A
\emph{warped} $G_2$-fundamental three-form $\phi$ is defined on $M$ by
\begin{equation}\label{eq:53}
  \phi=\omega_t\wedge dt + \psi_t^+.
\end{equation}
This is compatible with the warped product metric $g = dt^2+f^2g^*$
and has $ \Hodge\phi = \tfrac12\omega_t^2 + \psi_t^-\wedge dt$.

Suppose $(M^*,g^*,J,\omega,\psi^+,\psi^-)$ is nearly K\"ahler,
normalized according to equation~\eqref{eq:51} with structure
equations~\eqref{eq:52}.  We define the warped $G_2$ structure as in
equation~\eqref{eq:53}.  Then
\begin{gather*}
  d\omega_t = 2f^{-1}f'dt\wedge\omega_t +
  3f^{-1}\sigma(\cos\theta\psi^+_t+\sin\theta\psi_t^-),\\
  d\psi^+_t = 3f^{-1}f'dt\wedge\psi^+_t - d\theta\wedge\psi^-_t +
  2f^{-1}\sigma\sin\theta\omega_t^2,\\
  d\psi^-_t = 3f^{-1}f'dt\wedge\psi^-_t + d\theta\wedge\psi^+_t -
  2f^{-1}\sigma\cos\theta\omega_t^2,
\end{gather*}
and the differentials of the fundamental forms are
\begin{gather*}
  d\phi = -3f^{-1}\left(f' - \sigma\cos\theta\right)\psi^+_t\wedge dt
   + \left(\theta' + 3f^{-1}\sigma\sin\theta\right)\psi^-_t\wedge
  dt+ 2f^{-1}\sigma\sin\theta\omega_t^2,\\
  d\Hodge\phi = 2f^{-1}\left(f' -
  \sigma\cos\theta\right)dt\wedge\omega_t^2. 
\end{gather*}
From this the torsion components $\tau_p\in\Omega^p(M)$ are obtained.
We have
\begin{equation*}
  \begin{array}{ll}
    \tau_0 = \tfrac47\left(\theta' + 6f^{-1}\sigma\sin\theta
    \right),&
    \tau_1 = f^{-1}\left(f' - \sigma\cos\theta\right) dt,\\
    \tau_2 = 0,&
    \tau_3 = -\tfrac17\left(\theta' - f^{-1}\sigma\sin\theta\right)
    \left(4\omega_t\bw dt - 3\psi_t^+\right).
  \end{array}
\end{equation*}
The following elementary fact is recorded here for ease of reference
\begin{lemma}\label{lem:ode1}
  Suppose $b$ is a non-zero continuous function.  Then the solutions
  to the equation $\theta'=b\sin\theta$ are:
  either $\sin\theta=0=\theta'$, or
  $  \cos\theta=\frac{1-a^2}{1+a^2},\quad \sin\theta=\pm\frac{2a}{1+a^2},
  $
  where $a(t)=\exp{\int^t b(s) ds}$.
\end{lemma}
This Lemma ensures that given any function $f$, the torsion components
$\tau_0$ and $\tau_3$ may be made to vanish either simultaneously,
with $\sin\theta=0=\theta'$ or, when $\sigma\not=0$, separately, by
choosing an appropriate solution $\theta$ to the equation
$\theta'=b\sin\theta$ above.

\begin{proposition}\label{prop:6}
  Suppose $(M^*,g^*,\omega,\psi^+,\psi^-)$ is a nearly K\"ahler
  manifold.  Then any warped product $M=I\times M^*$, $g=dt^2+f^2 g^*$
  has $\Ric^\mW = 0$.  Any Einstein warped product $M=I\times M^*$,
  $g=dt^2+f^2 g^*$ has $\Ric^g_0 = 0 = \Ric^\phi_0$ and $W_{64}=0$.
\end{proposition}
\begin{proof}
  The first statement follows by using the formulas for the components
  $\tau_p$ of the torsion of $\phi$ computed above in the expression
  for the Ricci curvature given by Lemma~\ref{prop:1}.  The second is
  an easy consequence of the forms of the curvature tensors of nearly
  K\"ahler manifolds and Einstein warped product metrics.
\end{proof}

\subsection{Cohomogeneity one examples}

A variation on this theme is obtained by considering a homogeneous
$M^*=G/H$ with $H\subset\SU(3)$.  To be explicit we take
$G=\SU(3)$,~$H=T^2$.  On $M^*$ there are invariant forms
$\omega_i\in\Lambda^2,~i=1,2,3$, $\psi^+,\psi^-\in\Lambda^3$ such that
\begin{gather*}
  \omega_i\bw\psi^\pm=0,\quad \omega_i^2=0,\quad
  \vol_0:=\omega_1\omega_2\omega_3=\tfrac14\psi^+\bw\psi^-,\quad
  d\omega_i=\tfrac12\psi^+,\qquad d\psi^-=-2\sum_{i<j}\omega_i\omega_j.
\end{gather*}
On $M=I\times M^*$ we set
\begin{gather*}
  \omega_{i,t}:=f_i^2\omega_i,\quad
  \psi^+_t:=f_1f_2f_3(\cos\theta\psi^+-\sin\theta\psi_-),\quad
  \psi^-_t:=f_1f_2f_3(\sin\theta\psi^++\cos\theta\psi_-),\\
  \omega_t:=\sum_i\omega_{i,t},\quad
  \vol_t=\tfrac14\psi^+_t\bw\psi^-_t=\tfrac16\omega^3_t,\quad
  \Omega_t:=f_1f_2f_3\sum_{i<j}\omega_i\omega_j =
  \sum_{i<j}\tfrac{f_k}{f_if_j}\omega_{i,t}\omega_{j,t},
\end{gather*}
with $\{i,j,k\}=\{1,2,3\}$.  An easy calculation gives
\begin{gather*}
  d\omega_{i,t}=\ln(f_i^2)'dt\bw\omega_{i,t} +
  \tfrac12\tfrac{f_i}{f_jf_k}(\cos\theta\psi^+_t+\sin\theta\psi^-_t),\\
  d\omega_t = \sum_i\ln(f_i^2)'dt\bw\omega_{i,t} +
  h(\cos\theta\psi^+_t+\sin\theta\psi^-_t),\\
  d\psi^+_t = -k'\psi^+_t \bw dt + \theta'\psi^-_t\bw dt +
  2\sin\theta\Omega_t,\quad
  d\psi^-_t = -k'\psi^-_t \bw dt - \theta'\psi^+_t\bw dt
  - 2\cos\theta\Omega_t, 
\end{gather*}
where $k=\ln(f_1f_2f_3),~h=\tfrac{f_1^2+f_2^2+f_3^2}{2f_1f_2f_3}$. Now
set $\phi:=\omega_t\bw dt + \psi^+_t$.  Then $\phi$ is compatible with
the metric $g=dt^2 + \sum_if_i^2 g_i,$ and volume $\vol=\vol_t\bw dt,
$\hspace{2mm} and $\Hodge\phi =\tfrac12\omega_t^2+\psi_t^-\bw dt$.

Suppose that $f_1,f_2,f_3$ are such that $(f_if_j)'=f_k$ (there is a
one-parameter family of such triples, see~\cite{MR1969782} for
details).  Then the metric $g$ has holonomy contained in $G_2$.  For
such a triple one furthermore has $k'=h$.  Taking this in to account we
calculate the torsion components of $\phi$.
\begin{gather*}
  \tau_0 := \tfrac47(\theta'+2h\sin\theta),\\
  \tau_1 := \tfrac13h(1-\cos\theta),\\
  \tau_2 := -\frac{2(1 - \cos\theta)}{3f_1f_2f_3}\sum_i(2f_i^2 - f_j^2
  - f_k^2) \omega_{i,t},\\
  \tau_3:=\tfrac37\left(\theta' - \tfrac13h\sin\theta\right)\psi_t^+ -
  \tfrac47\sum_i\left(\theta' - \frac{5f^2_i - 2(f_j^2 +
      f_k^2)}{2f_1f_2f_3}\sin\theta\right)\omega_{i,t}\bw dt.
\end{gather*}
First note that $\phi$ is parallel if and only if $\cos\theta=1$.  If
$\cos\theta\not\equiv 1$ then $\tau_1\not\equiv0$.  Taking
$\cos\theta=-1$ gives a $2+4$ structure, which is a type $4$ structure
when $f_1=f_2=f_3$.

If $f_1=f_2=f_3=:f$ holds then $\tau_2\equiv0$ and \[ \tau_3=
\left(\theta' - \tfrac13h\sin\theta\right)\left(\tfrac37\psi_t^+ -
  \tfrac47\omega_t\bw dt\right).\] In general this gives a type
$1+3+4$ structure, but using Lemma~\ref{lem:ode1} again we may make
either $\tau_0$ or $\tau_3$ vanish to obtain $1+4$ and $3+4$
structures, too.

Now suppose $f_3\geqslant f_2\geqslant f_1$ and $f_3 > f_1$.  The
generic type is $1+2+3+4$, but we may once again use
Lemma~\ref{lem:ode1} to eliminate $\tau_0$ to get class $2+3+4$.
Proposition~\ref{prop:6} and a little book-keeping then proves
\begin{proposition}\label{prop:8}
  There are warped product and cohomogeneity-one metrics $g$ with
  $\Hol(g)\subset G_2$ and three-forms $\phi$ compatible with $g$ of
  every admissible type, except possibly for $1+2+3$.\qed
\end{proposition}
A three-form of \emph{admissible type} means not one of the strict
types $1,2,3,$ and $2+3$.  Type $1+2$ doesn't really exist as a strict
class and every other type of $G_2$ structure in this list has either
strictly positive or strictly negative scalar curvature,
c.f.~\eqref{eq:29}.
\begin{remark}
  While the Riemannian curvature is a function of the metric alone,
  the $G_2$ equivariant splitting of the curvature is inherently
  dependent on the three-form, in general.  However, for the very
  special case of a warped product three-form with Einstein metric on
  the product of a nearly K\"ahler manifold with $\mathbb R$, the
  component $W_{77}$ equals the curvature component $K^*$ of the
  nearly K\"ahler manifold for any $\phi$.  For the cohomogeneity-one
  $G_2$ structures the picture is far more complicated as the constant
  $t$ hypersurfaces $(N^*,\sum_i f_i^2g_i)$ are half-flat rather than
  nearly K\"ahler.
\end{remark}
\begin{remark} The following was explained to us by Bryant
  \cite{priv:Br}: for a fixed metric $g$, the exterior differential
  system corresponding to the equation $p^4_7(d\phi)=0$ for a
  compatible $G_2$-three-form $\phi$ is involutive at points where the
  torsion $\tau\not\equiv0$, with last nonzero Cartan character $s_6 =
  6$, so that the general local solution depends on 6 functions of 6
  variables.  This along with Proposition~\ref{prop:8} entails
  \begin{theorem}
    For every admissible Fern\'andez-Gray type, there is a $G_2$
    three-form $\phi$, with torsion of this type such that the
    associated metric $g$ has holonomy $\Hol(g)\subset G_2$.
  \end{theorem}
\end{remark}

\subsection{Complete and compact examples}
\begin{example}\label{ex:1}
  Suppose $(M^*,\omega,\psi^+)$ is Calabi-Yau (corresponding to
  structure equations $d\omega=0, ~d\psi^+=0=d\psi^-$).  Then for any
  smooth function $\theta\colon\mathbb R\to\mathbb R$, the three-form
  $\phi$ defined as in \eqref{eq:53} with $f$ constant is smooth.  The $G_2$
  three-form is of strict type $1+3$ if $\theta$ is non-constant and
  parallel otherwise.  The associated metric is the Riemannian product
  metric on $\mathbb R\times M^*$ and hence has holonomy contained in
  $\SU(3)\subset G_2$.  Compact examples may be obtained by the method
  of~\cite{Besse:Einstein},~Section 9.109.
\end{example}
For certain other choices of $(M^*,g^*)$ and non-constant function
$f\colon I\to \mathbb R$ the warped product metric on $I\times M^*$
also extends to a complete metric.  This is the case when
\begin{align*}
  g_{\mathbb R^7} &= dt^2+t^2g_{S^6},\qquad\qquad I=\mathbb R^+,\\
  g_{S^7}&=dt^2+\sin^2(t)g_{S^6},\qquad I=(0,\pi),\\
  g_{H^7}&=
  \begin{cases}
    dt^2+\sinh^2(t)g_{S^6},& I=\mathbb R^+,\\
    dt^2+e^{2t}g_{\mathbb R^6},& I=\mathbb R.
  \end{cases}
\end{align*}
Here $g_N$ refers to the constant sectional curvature metric of $N$
with $\abs{\rho}$ and $\rho^*$ equal to $0$ or $1$.

\begin{example}
  We consider hyperbolic space $H^7$ with metric
  $dt^2+e^{2t}g_{\mathbb R^6}$ first.  The structure on $M^*$, as the
  examples of~\ref{ex:1}, has $\sigma=0$ and $f$ strictly positive.
  The warped three-form $\phi$ is therefore smooth for any smooth
  function $\theta\colon\mathbb R\to\mathbb R$.  The generic type is
  $1+3+4$ for non-constant $\theta$ and $4$ for constant $\theta$.
  Type $1+3$ three-forms require $e^t=\cos\theta$ and therefore are
  not complete.
\end{example}

The remaining three cases are all cohomogeneity-one metrics on the
isotropy irreducible space $M^*=S^6=G_2/\SU(3)$ equipped with its
standard, homogeneous, nearly K\"ahler structure so that $\sigma=1$.
The structures with $\tau_1=0$ are: the flat $G_2$ structure on
$\mathbb R^7$, where $f=t,~\theta=0$ and the unique nearly parallel
(type $1$) structure on $S^7$, $f=\sin(t),~\theta=t$.  These are, of
course, complete and smooth.

Type $1+4$ structures are given by the equation $fd\theta= \sin\theta
dt$.  Consequently, the metric and three-form may conveniently be
re-written in the familiar form
\begin{gather*}
  \phi=(\theta')^{-3}(\sin^2\theta\omega\bw d\theta +
  \sin^3\theta(\cos\theta\psi^+-\sin\theta\psi^-)),\qquad
  g=(\theta')^{-2}(d\theta^2+\sin^2\theta g_{S^6}).
\end{gather*}
It is then clear that the three-forms of type $1+4$ arise as
point-wise conformal changes of the standard nearly parallel structure
on the $7$-sphere.  Examination of the solutions of
Lemma~\ref{lem:ode1} corresponding to $b=f^{-1}$ shows that there are
constants $A$ and $B$ so that
\begin{equation*}
  \theta'=A\cos\theta + B,\quad\text{where}\quad
  \begin{cases}
    \abs{B/A}>1,&\text{for $f=\sin(t)$},\\
    \abs{B/A}<1,&\text{for $f=\sinh(t)$},\\
    \abs{B/A}=1,&\text{for $f=t$}.
  \end{cases}
\end{equation*}
The requirement that $\abs{\theta'}>0$ imposes restrictions in the two
latter cases.  It is clear that both metric and three-form are smooth
where $\theta'\not=0$.  In the first case we see that all type $1+4$
warped $G_2$ structures arise as global conformal changes of the
standard nearly parallel structure on $S^7$,
see~\cite{math.DG/0607487,MR1512246}.

\begin{example}
  Type $4$ structures are also interesting. First consider three-forms
  $\phi$ compatible with the standard metric on $\mathbb R^7$.  Since
  $\theta'=0=\sin\theta$, there are only two possibilities: either
  $\phi$ is parallel ($\cos\theta=1$) or $\tau_1 = d\ln(t^2)$
  ($\cos\theta=-1$).  In the latter case one notes that
  \begin{gather*}
    t^{-6}\phi = s^2\omega\bw ds + s^3\psi^+ = \iota^*\phi,\qquad
    t^{-4}g =ds^2+s^2 g_{S^6} = \iota^* g,
  \end{gather*}
  where $s:=-t^{-1}$ and $\iota$ is the map $\iota\colon\mathbb
  R^7\setminus\{0\}\to \mathbb R^7\setminus\{0\}, ~x \mapsto
  -x/\abs{x}^2$.

  The standard metric of the $7$ sphere is compatible with the
  three-form
  \begin{gather*}
    \phi=\sin^2(t)\omega\bw dt + \sin^3(t)\psi^+.
  \end{gather*}
  We note that this three-form satisfies the necessary conditions
  of~\cite{MR1969782}.  However, the three-form has type $4$ with
  torsion $\tau_1=-\tan(t/2)dt=-d\ln\cos^2(t/2)$.  This is clearly
  singular at one point.  Setting $r=2\cos^2(t/2),~s=\tan(t/2)$ then
  $ds=\tfrac12 dt/\cos^2(t/2)$ and we may write
  \begin{gather*}
    \phi=\frac8{(1+s^2)^3}(s^2\omega\bw ds + s^3\psi^+),\qquad
    g_{S^7}=\frac4{(1+s^2)^2}(ds^2+s^2g_{S^6}) = r^2g_{\mathbb R^7},
  \end{gather*}
  This transformation realizes the same structure smoothly on $\mathbb
  R^7$. The metric represented this way is not complete.
  
  For the hyperbolic space the same technique produces
  \begin{gather*}
    \phi=\frac8{(1-s^2)^3}(s^2\omega\bw ds + s^3\psi^+),\qquad
    g=\frac4{(1-s^2)^2}(ds^2+s^2 g_{S^6}),
  \end{gather*}
  where $s=\tanh(t/2)\in(0,1)$.  
\end{example}

\subsubsection{Cohomogeneity one metrics over $\SU(3)/T^2$}

The analysis carried out in the final section of~\cite{MR1969782}
applies to the cohomogeneity one metrics based on
$M^*=\SU(3)/T^2,\SP(2)/\SP(1)\Un(1)$, too.  Whenever $\theta$ is an
odd function of $t$, $f_1$ is odd with $f_1'(0)=1$ and $f_2$ is even
with $f_2(0)>0$ the metric $g$ and three-form $\phi$ extend smoothly
to the non-compact manifolds isomorphic to the bundles of
anti-self-dual forms over $\mathbb CP(2)$ and $S^4$, respectively.
These then carry smooth $G_2$ structures compatible with the holonomy
$G_2$ structures which are either parallel or of strict type $2+4$,
$2+3+4$, $1+2+3+4$.

\def\cprime{$'$}
\providecommand{\bysame}{\leavevmode\hbox to3em{\hrulefill}\thinspace}
\providecommand{\href}[2]{#2}

\end{document}